\documentclass[11pt,twoside]{article}

\usepackage[latin1]{inputenc}
\usepackage{epsfig}
\usepackage{amsfonts}
\usepackage{cite}
\usepackage{color}
\definecolor{lgrey}{gray}{0.99}

\title{
  {\huge Complex Analysis of Real Functions \\[1.5ex]}
  VI: On the Convergence of Fourier Series }

\author{
  \Large Jorge L. deLyra\footnote{Email: delyra@latt.if.usp.br} \\
  Department of Mathematical Physics \\
  Physics Institute \\
  University of São Paulo }

\date{April 8, 2018}

\newcommand{\ii}{\mbox{\boldmath$\imath$}}
\newcommand{\iii}{\mbox{\boldmath\scriptsize$\imath$}}

\newcommand{\HH}{{\cal H}_{\rm c}}
\newcommand{\DS}{{\cal D}_{\rm s}}
\newcommand{\DC}{{\cal D}_{\rm c}}
\newcommand{\DR}{{\cal D}_{\rm r}}

\newcommand{\e}[1]{\,{\rm e}^{#1}}

\newcommand{\ldot}{\mbox{\Large$\cdot$}\!}

\newcommand{\PV}{\mbox{\rm PV}\!\!}

\newtheorem{definition}{Definition}

\newtheorem{theorem}{Theorem}

\newtheorem{proof}{Proof}[theorem]
\newcommand{\Colon}{{\hspace{-0.4em}\bf:}}

\begin{document}\maketitle

\vspace{-1.65ex}
\begin{abstract}
  \noindent
  We define a compact version of the Hilbert transform, which we then use
  to write explicit expressions for the partial sums and remainders of
  arbitrary Fourier series. The expression for the partial sums reproduces
  the known result in terms of Dirichlet integrals. The expression for the
  remainder is written in terms of a similar type of integral. Since the
  asymptotic limit of the remainder being zero is a necessary and
  sufficient condition for the convergence of the series, this same
  condition on the asymptotic behavior of the corresponding integrals
  constitutes such a necessary and sufficient condition.
\end{abstract}

\section{Introduction}\label{Sec01}

In a previous paper~\cite{CAoRFI} we introduced a certain complex-analytic
structure within the unit disk of the complex plane, and showed that it is
possible to represent within that structure essentially {\em all\/}
integrable real functions defined in a compact interval. In a subsequent
paper~\cite{CAoRFIII} we showed that {\em all\/} the elements of the
Fourier theory~\cite{FSchurchill} of integrable real functions are
contained within that complex-analytic structure. However, in that paper
we did not discuss in any depth the question of the convergence of Fourier
series.

The fact that it is possible to recover the real functions from their
Fourier coefficients almost everywhere, even when the corresponding
Fourier series are divergent, as we showed in~\cite{CAoRFI}, led to a
powerful and very general summation rule for {\em all\/} Fourier series,
which was presented in~\cite{CAoRFIII}. This summation rule allows one to
add up a regularized version of the Fourier series, in a meaningful way,
and therefore allows one to simply circumvent the fact that the original
Fourier series may be divergent. However, the complex-analytic structure
actually does allow for a direct discussion of the convergence problem.

In this paper we will present a more complete analysis of the convergence
of Fourier series. In order to do this we will first introduce what we
will name the {\em compact Hilbert transform}, which is a version of the
Hilbert transform which is appropriate for functions defined on a compact
interval. This will lead not only to the known explicit expression for the
partial sums of the Fourier series in terms of Dirichlet integrals, but
also to an explicit expression for the {\em remainder} of the Fourier
series, in terms of a similar type of integral.

For ease of reference, we include here a one-page synopsis of the
complex-analytic structure introduced in~\cite{CAoRFI}. It consists of
certain elements within complex analysis~\cite{CVchurchill}, as well as of
their main properties.

\paragraph{Synopsis:} The Complex-Analytic Structure\\

\noindent
An {\em inner analytic function} $w(z)$ is simply a complex function which
is analytic within the open unit disk. An inner analytic function that has
the additional property that $w(0)=0$ is a {\em proper inner analytic
  function}. The {\em angular derivative} of an inner analytic function is
defined by

\noindent
\begin{equation}
  w^{\ldot}(z)
  =
  \ii
  z\,
  \frac{dw(z)}{dz}.
\end{equation}

\noindent
By construction we have that $w^{\ldot}(0)=0$, for all $w(z)$. The {\em
  angular primitive} of an inner analytic function is defined by

\begin{equation}
  w^{-1\ldot}(z)
  =
  -\ii
  \int_{0}^{z}dz'\,
  \frac{w(z')-w(0)}{z'}.
\end{equation}

\noindent
By construction we have that $w^{-1\ldot}(0)=0$, for all $w(z)$. In terms
of a system of polar coordinates $(\rho,\theta)$ on the complex plane,
these two analytic operations are equivalent to differentiation and
integration with respect to $\theta$, taken at constant $\rho$. These two
operations stay within the space of inner analytic functions, they also
stay within the space of proper inner analytic functions, and they are the
inverses of one another. Using these operations, and starting from any
proper inner analytic function $w^{0\ldot}(z)$, one constructs an infinite
{\em integral-differential chain} of proper inner analytic functions,

\begin{equation}
  \left\{
    \ldots,
    w^{-3\ldot}(z),
    w^{-2\ldot}(z),
    w^{-1\ldot}(z),
    w^{0\ldot}(z),
    w^{1\ldot}(z),
    w^{2\ldot}(z),
    w^{3\ldot}(z),
    \ldots\;
  \right\}.
\end{equation}

\noindent
Two different such integral-differential chains cannot ever intersect each
other. There is a {\em single} integral-differential chain of proper inner
analytic functions which is a constant chain, namely the null chain, in
which all members are the null function $w(z)\equiv 0$.

A general scheme for the classification of all possible singularities of
inner analytic functions is established. A singularity of an inner
analytic function $w(z)$ at a point $z_{1}$ on the unit circle is a {\em
  soft singularity} if the limit of $w(z)$ to that point exists and is
finite. Otherwise, it is a {\em hard singularity}. Angular integration
takes soft singularities to other soft singularities, and angular
differentiation takes hard singularities to other hard singularities.

Gradations of softness and hardness are then established. A hard
singularity that becomes a soft one by means of a single angular
integration is a {\em borderline hard} singularity, with degree of
hardness zero. The {\em degree of softness} of a soft singularity is the
number of angular differentiations that result in a borderline hard
singularity, and the {\em degree of hardness} of a hard singularity is the
number of angular integrations that result in a borderline hard
singularity. Singularities which are either soft or borderline hard are
integrable ones. Hard singularities which are not borderline hard are
non-integrable ones.

Given an integrable real function $f(\theta)$ on the unit circle, one can
construct from it a unique corresponding inner analytic function $w(z)$.
The real function $f(\theta)$ is recovered by means of the $\rho\to
1_{(-)}$ limit of the real part of this inner analytic function.
Singularities of real functions can be classified in a way which is
analogous to the corresponding complex classification. Integrable real
functions are typically associated with inner analytic functions that have
singularities which are either soft or at most borderline hard. A more
detailed review of real functions will be given in
Section~\ref{Sec02}. This ends our synopsis.

\vspace{2.6ex}

\noindent
Some of the material contained in this paper can be seen as a development,
reorganization and extension of some of the material found, sometimes
still in rather rudimentary form, in the
papers~\cite{FTotCPI,FTotCPII,FTotCPIII,FTotCPIV,FTotCPV}.

\section{Review of Real Functions}\label{Sec02}

When we discuss real functions in this paper, some properties will be
globally assumed for these functions, just as was done
in~\cite{CAoRFI,CAoRFII,CAoRFIII}. These are rather weak conditions to be
imposed on these functions, that will be in force throughout this paper.
It is to be understood, without any need for further comment, that these
conditions are valid whenever real functions appear in the arguments.
These weak conditions certainly hold for any integrable real functions
that are obtained as restrictions of corresponding inner analytic
functions to the unit circle.

The most basic condition is that the real functions must be measurable in
the sense of Lebesgue, with the usual Lebesgue
measure~\cite{RARudin,RARoyden}. The second global condition we will
impose is that the functions have no removable singularities. The third
and last global condition is that the number of hard singularities on the
unit circle be finite, and hence that they be all isolated from one
another. There will be no limitation on the number of soft singularities.

In addition to this we will assume, for the purposes of this particular
paper, that all real functions are integrable on the unit circle and, just
for the sake of clarity and simplicity, unless explicitly stated otherwise
we will also assume that all real functions are zero-average real
functions, meaning that their integrals over the unit circle are zero.
Since this simply implies that the Fourier coefficients $\alpha_{0}$ of
the real functions are zero, without affecting any of the other
coefficients in any way, this clearly has no impact on any arguments about
the convergence of the Fourier series.

For the purposes of this paper it is important to review here, in some
detail, the construction that results in the correspondence between
integrable real functions on the unit circle and inner analytic functions
on the open unit disk. In~\cite{CAoRFI} we showed that, given any
integrable real function $f(\theta)$, one can construct a corresponding
inner analytic function $w(z)$, from the real part of which $f(\theta)$
can be recovered almost everywhere on the unit circle, through the use of
the $\rho\to 1_{(-)}$ limit, where $(\rho,\theta)$ are polar coordinates
on the complex plane. In that construction we started by calculating the
Fourier coefficients~\cite{FSchurchill} $\alpha_{k}$ and $\beta_{k}$ of
the real function, which is always possible given that the function is
integrable, using the usual integrals defining these coefficients,

\noindent
\begin{eqnarray}\label{EQFourCoef1}
  \alpha_{0}
  & = &
  \frac{1}{\pi}
  \int_{-\pi}^{\pi}d\theta\,
  f(\theta),
  \nonumber\\
  \alpha_{k}
  & = &
  \frac{1}{\pi}
  \int_{-\pi}^{\pi}d\theta\,
  \cos(k\theta)f(\theta),
  \nonumber\\
  \beta_{k}
  & = &
  \frac{1}{\pi}
  \int_{-\pi}^{\pi}d\theta\,
  \sin(k\theta)f(\theta),
\end{eqnarray}

\noindent
for $k\in\{1,2,3,\ldots,\infty\}$. We then defined a set of complex Taylor
coefficients $c_{k}$ given by

\noindent
\begin{eqnarray}\label{EQTaylCoef}
  c_{0}
  & = &
  \frac{1}{2}\,\alpha_{0},
  \nonumber\\
  c_{k}
  & = &
  \alpha_{k}
  -
  \ii\beta_{k},
\end{eqnarray}

\noindent
for $k\in\{1,2,3,\ldots,\infty\}$. Next we defined a complex variable $z$
associated to $\theta$, using the positive real variable $\rho$, by
$z=\rho\exp(\ii\theta)$. Using all these elements we then constructed the
complex power series

\begin{equation}\label{EQPowrSers}
  S(z)
  =
  \sum_{k=0}^{\infty}
  c_{k}z^{k},
\end{equation}

\noindent
which we showed in~\cite{CAoRFI} to be convergent to an inner analytic
function $w(z)$ within the open unit disk. That inner analytic function
may be written as

\begin{equation}\label{EQInneAnalFunc}
  w(z)
  =
  u(\rho,\theta)
  +
  \ii
  v(\rho,\theta).
\end{equation}

\noindent
The complex power series in Equation~(\ref{EQPowrSers}) is therefore the
Taylor series of $w(z)$. We also proved in~\cite{CAoRFI} that one recovers
the real function $f(\theta)$ almost everywhere on the unit circle from
the real part $u(\rho,\theta)$ of $w(z)$, by means of the $\rho\to
1_{(-)}$ limit. The $\rho\to 1_{(-)}$ limit of the imaginary part
$v(\rho,\theta)$ also exists almost everywhere and gives rise to a real
function $g(\theta)$ which corresponds to $f(\theta)$. The pair of real
functions obtained from the real and imaginary parts of one and the same
inner analytic function are said to be mutually Fourier-conjugate real
functions.

In a subsequent paper~\cite{CAoRFIII} we showed that {\em all\/} the
elements of the Fourier theory~\cite{FSchurchill} of integrable real
functions are contained within the complex-analytic structure, including
the Fourier basis of functions, the Fourier series, the scalar product for
integrable real functions, the relations of orthogonality and norm of the
basis elements, and the completeness of the Fourier basis, including its
so-called completeness relation. As was also shown in~\cite{CAoRFIII} the
real function $g(\theta)=v(1,\theta)$ which is the Fourier-conjugate
function to $f(\theta)=u(1,\theta)$ has the same Fourier coefficients, but
with the meanings of $\alpha_{k}$ and $\beta_{k}$ interchanged in such a
way that we have

\noindent
\begin{eqnarray}\label{EQFourCoef2}
  \alpha_{k}
  & = &
  \frac{1}{\pi}
  \int_{-\pi}^{\pi}d\theta\,
  \sin(k\theta)
  g(\theta),
  \nonumber\\
  \beta_{k}
  & = &
  -\,
  \frac{1}{\pi}
  \int_{-\pi}^{\pi}d\theta\,
  \cos(k\theta)
  g(\theta),
\end{eqnarray}

\noindent
for $k\in\{1,2,3,\ldots,\infty\}$. Note that there is no constant term in
the Fourier series of $g(\theta)$, which means that we have

\begin{equation}
  \int_{-\pi}^{\pi}d\theta\,
  g(\theta)
  =
  0.
\end{equation}

\noindent
In other words, the Fourier-conjugate function $g(\theta)$ is always a
zero-average real function. Note also that this fact, as well as the
relations in Equation~(\ref{EQFourCoef2}) imply, in particular, that
$g(\theta)$ is also an integrable real function. We may therefore conclude
that, if $f(\theta)$ is an integrable real function, then so is its
Fourier-conjugate function $g(\theta)$.

\section{The Compact Hilbert Transform}\label{Sec03}

Let $f(\theta)$ be an integrable real function on $[-\pi,\pi]$, with
Fourier coefficients as given in Equation~(\ref{EQFourCoef1}). As was
shown in~\cite{CAoRFIII} the real function $g(\theta)=v(1,\theta)$ which
is the Fourier-conjugate function to $f(\theta)=u(1,\theta)$ has the same
Fourier coefficients, but with the meanings of $\alpha_{k}$ and
$\beta_{k}$ interchanged, as shown in Equation~(\ref{EQFourCoef2}). The
relations in Equations~(\ref{EQFourCoef1}) and~(\ref{EQFourCoef2}) can be
understood as the following collection of integral identities satisfied by
all pairs of Fourier-conjugate integrable real functions,

\noindent
\begin{eqnarray}
  \int_{-\pi}^{\pi}d\theta\,
  \cos(k\theta)
  f(\theta)
  & = &
  \int_{-\pi}^{\pi}d\theta\,
  \sin(k\theta)
  g(\theta),
  \nonumber\\
  \int_{-\pi}^{\pi}d\theta\,
  \sin(k\theta)
  f(\theta)
  & = &
  -\,
  \int_{-\pi}^{\pi}d\theta\,
  \cos(k\theta)
  g(\theta),
\end{eqnarray}

\noindent
for $k\in\{1,2,3,\ldots,\infty\}$. It is well known that this replacement
of $\cos(k\theta)$ with $\sin(k\theta)$ and of $\sin(k\theta)$ with
$-\cos(k\theta)$ can be effected by the use of the Hilbert transform.
However, that transform was originally introduced by Hilbert for real
functions defined on the whole real line, rather that on the unit circle
as is our case here. Therefore, the first thing that we will do here is to
define a compact version of the Hilbert transform that applies to real
functions defined on the unit circle.

Since the Fourier coefficient $\alpha_{0}$ of $f(\theta)$ has no effect on
the definition of the Fourier-conjugate function $g(\theta)$, and in order
for this pair of real functions to be related in a unique way, we will
assume that $f(\theta)$ is also a zero-average real function,

\begin{equation}
  \int_{-\pi}^{\pi}d\theta\,
  f(\theta)
  =
  0,
\end{equation}

\noindent
thus implying for its $k=0$ Fourier coefficient that $\alpha_{0}=0$. This
does not affect, of course, any subsequent arguments about the convergence
of the Fourier series. According to the construction presented
in~\cite{CAoRFI} and reviewed in Section~\ref{Sec02}, from the other
Fourier coefficients $\alpha_{k}$ and $\beta_{k}$, for
$k\in\{1,2,3,\ldots,\infty\}$, we may construct the complex coefficients
$c_{k}=\alpha_{k}-\ii\beta_{k}$, for $k\in\{1,2,3,\ldots,\infty\}$, where
we now have $c_{0}=0$, and from these we may construct the corresponding
inner analytic function $w(z)$ shown in Equation~(\ref{EQInneAnalFunc}),
which is now, in fact, a {\em proper} inner analytic function, since
$c_{0}=0$ implies that $w(0)=0$.

Since $u(\rho,\theta)$ and $v(\rho,\theta)$ are harmonic conjugate
functions to each other, it is now clear that there is a one-to-one
correspondence between $u(\rho,\theta)$ and $v(\rho,\theta)$, and in
particular between $u(1,\theta)$ and $v(1,\theta)$. Therefore, there is a
one-to-one correspondence between $f(\theta)$ and $g(\theta)$, in this
case valid almost everywhere on the unit circle, since we have shown
in~\cite{CAoRFI} that $f(\theta)=u(1,\theta)$ and that
$g(\theta)=v(1,\theta)$, both almost everywhere over the unit circle.
Therefore, a transformation must exist that produces $g(\theta)$ from
$f(\theta)$ almost everywhere over the unit circle, as well as an inverse
transformation that recovers $f(\theta)$ from $g(\theta)$ almost
everywhere over the unit circle. In this section we will show that the
following definition accomplishes this.

\begin{definition}\Colon\label{Def01}
  Compact Hilbert Transform
\end{definition}

\noindent
If $f(\theta)$ is an arbitrarily given zero-average integrable real
function defined on the unit circle, then its {\em compact Hilbert
  transform} $g(\theta)$ is the real function defined by

\noindent
\begin{eqnarray}
  g(\theta)
  & = &
  \HH[f(\theta)]
  \nonumber\\
  & = &
  -\,
  \frac{1}{2\pi}\,
  \PV\int_{-\pi}^{\pi}d\theta_{1}\,
  \frac
  {\cos\!\left[\rule{0em}{2ex}(\theta_{1}-\theta)/2\right]}
  {\sin\!\left[\rule{0em}{2ex}(\theta_{1}-\theta)/2\right]}\,
  f(\theta_{1}),
\end{eqnarray}

\noindent
where $\PV\,$ stands for the Cauchy principal value, and where
$\HH[f(\theta)]$ is the notation we will use for the compact Hilbert
transform applied to the real function $f(\theta)$.

\vspace{2.6ex}

\noindent
We will now prove the following theorem.

\begin{theorem}\Colon\label{Theo01}
  The zero-average integrable real functions $f(\theta)$ and $g(\theta)$,
  which are such that $f(\theta)=u(1,\theta)$ and $g(\theta)=v(1,\theta)$
  almost everywhere on the unit circle, are related to each other by this
  transform, that is, we have that $g(\theta)=\HH[f(\theta)]$ almost
  everywhere on the unit circle, and that $f(\theta)=\HH^{-1}[g(\theta)]$
  almost everywhere on the unit circle, where the inverse transform is
  simply given by $\HH^{-1}[g(\theta)]=-\HH[g(\theta)]$.
\end{theorem}

\newpage

\begin{proof}\Colon
\end{proof}

\noindent
In order to derive these facts from our complex-analytic structure, we
start from the Cauchy integral formula for the inner analytic function
$w(z)$,

\begin{equation}\label{EQCauIntForI}
  w(z)
  =
  \frac{1}{2\pi\ii}
  \oint_{C}dz_{1}\,
  \frac{w(z_{1})}{z_{1}-z},
\end{equation}

\noindent
where $C$ can be taken as a circle centered at the origin, with radius
$\rho_{1}<1$, and where we write $z$ and $z_{1}$ in polar coordinates as
$z=\rho\exp(\ii\theta)$ and $z_{1}=\rho_{1}\exp(\ii\theta_{1})$. The
integral formula in Equation~(\ref{EQCauIntForI}) is valid for
$\rho<\rho_{1}$, and in fact, by the Cauchy-Goursat theorem, the integral
is zero if $\rho>\rho_{1}$, since both $z$ and $z_{1}$ are within the open
unit disk, a region where $w(z)$ is analytic. We must now determine what
happens if $\rho=\rho_{1}$, that is, if $z$ is on the circle $C$ of radius
$\rho_{1}$. Note that in this case we may slightly deform the integration
contour $C$ in order to have it pass on one side or the other of the
simple pole of the integrand at $z_{1}=z$. If we use a deformed contour
$C_{\ominus}$ that {\em excludes} the pole from its interior, then we
have, instead of Equation~(\ref{EQCauIntForI}),

\begin{equation}\label{EQExclPole1}
  0
  =
  \frac{1}{2\pi\ii}
  \oint_{C_{\ominus}}dz_{1}\,
  \frac{w(z_{1})}{z_{1}-z},
\end{equation}

\noindent
due to the Cauchy-Goursat theorem, while if we use a deformed contour
$C_{\oplus}$ that {\em includes} the pole in its interior, then we have,
just as in Equation~(\ref{EQCauIntForI}),

\begin{equation}\label{EQInclPole1}
  w(z)
  =
  \frac{1}{2\pi\ii}
  \oint_{C_{\oplus}}dz_{1}\,
  \frac{w(z_{1})}{z_{1}-z}.
\end{equation}

\noindent
Since by the Sokhotskii-Plemelj theorem~\cite{sokhplem} the Cauchy
principal value of the integral over the circle $C$ is the arithmetic
average of these two integrals, in the limit where the deformations
vanish, a limit which does not really have to be considered in detail, so
long as the deformations do not cross any other singularities,

\begin{equation}
  \PV\oint_{C}dz_{1}\,
  \frac{w(z_{1})}{z_{1}-z}
  =
  \frac{1}{2}
  \oint_{C_{\ominus}}dz_{1}\,
  \frac{w(z_{1})}{z_{1}-z}
  +
  \frac{1}{2}
  \oint_{C_{\oplus}}dz_{1}\,
  \frac{w(z_{1})}{z_{1}-z},
\end{equation}

\noindent
adding Equations~(\ref{EQExclPole1}) and~(\ref{EQInclPole1}) we may
conclude that

\begin{equation}\label{EQCauchyPV}
  w(z)
  =
  \frac{1}{\pi\ii}\,
  \PV\oint_{C}dz_{1}\,
  \frac{w(z_{1})}{z_{1}-z},
\end{equation}

\noindent
where we now have $\rho=\rho_{1}$, that is, both $z_{1}$ and $z$ are on
the circle $C$ of radius $\rho_{1}$ within the open unit disk. This
formula can be understood as a special version of the Cauchy integral
formula, and will be used repeatedly in what follows. We may now write all
quantities in this equation in terms of the polar coordinates $\rho_{1}$,
$\theta_{1}$ and $\theta$,

\noindent
\begin{eqnarray}
  w(\rho_{1},\theta)
  & = &
  \frac{1}{\pi\ii}\,
  \PV\int_{-\pi}^{\pi}d\theta_{1}\,
  \ii\rho_{1}\e{\iii\theta_{1}}\,
  \frac
  {w(\rho_{1},\theta_{1})}
  {\rho_{1}\e{\iii\theta_{1}}-\rho_{1}\e{\iii\theta}}
  \nonumber\\
  & = &
  \frac{1}{\pi}\,
  \PV\int_{-\pi}^{\pi}d\theta_{1}\,
  \frac
  {u(\rho_{1},\theta_{1})+\ii v(\rho_{1},\theta_{1})}
  {1-\e{-\iii\Delta\theta}},
\end{eqnarray}

\noindent
where $\Delta\theta=\theta_{1}-\theta$. Note that, since by construction
the real and imaginary parts $u(\rho_{1},\theta)$ and $v(\rho_{1},\theta)$
of $w(\rho_{1},\theta)$ for $\rho_{1}=1$ are integrable real functions on
the unit circle, and since there are no other dependencies on $\rho_{1}$
in this expression, we may now take the $\rho_{1}\to 1_{(-)}$ limit of
this equation, in which the principal value acquires its usual real
meaning over the unit circle, that is, the meaning that the asymptotic
limits of the integral on either side of a non-integrable singularity must
be taken in the symmetric way. In that limit we have
 
\begin{equation}\label{EQCauIntForII}
  u(1,\theta)+\ii v(1,\theta)
  =
  \frac{1}{\pi}\,
  \PV\int_{-\pi}^{\pi}d\theta_{1}\,
  \frac
  {u(1,\theta_{1})+\ii v(1,\theta_{1})}
  {1-\e{-\iii\Delta\theta}}.
\end{equation}

\noindent
In order to identify separately the real and imaginary parts of this
equation, we must now rationalize the integrand of the integral shown. We
will use the fact that

\noindent
\begin{eqnarray}
  \frac{1}{1-\e{-\iii\Delta\theta}}
  & = &
  \frac
  {
    1
    -
    \e{\iii\Delta\theta}
  }
  {
    \left(
      \rule{0em}{2ex}
      1
      -
      \e{-\iii\Delta\theta}
    \right)
    \left(
      \rule{0em}{2ex}
      1
      -
      \e{\iii\Delta\theta}
    \right)
  }
  \nonumber\\
  & = &
  \frac
  {
    \left[
      \rule{0em}{2ex}
      1
      -
      \cos(\Delta\theta)
    \right]
    -
    \ii
    \sin(\Delta\theta)
  }
  {
    2-2\cos(\Delta\theta)
  }
  \nonumber\\
  & = &
  \frac{1}{2}
  \left[
    1
    -
    \ii\,
    \frac
    {
      \sin(\Delta\theta)
    }
    {
      1-\cos(\Delta\theta)
    }
  \right].
\end{eqnarray}

\noindent
Using the half-angle trigonometric identities we may also write this
result as

\noindent
\begin{eqnarray}\label{EQRatiForm1}
  \frac{1}{1-\e{-\iii\Delta\theta}}
  & = &
  \frac{1}{2}
  \left[
    1
    -
    \ii\,
    \frac{\cos(\Delta\theta/2)}{\sin(\Delta\theta/2)}
  \right]
  \\\label{EQRatiForm2}
  & = &
  -\,
  \frac{\ii}{2}\,
  \frac
  {\e{\iii\Delta\theta/2}}
  {\sin(\Delta\theta/2)}.
\end{eqnarray}

\noindent
Using the result shown in Equation~(\ref{EQRatiForm1}) back in
Equation~(\ref{EQCauIntForII}) we obtain

\noindent
\begin{eqnarray}
  u(1,\theta)
  +
  \ii
  v(1,\theta)
  & = &
  \frac{1}{2\pi}\,
  \PV\int_{-\pi}^{\pi}d\theta_{1}\,
  \left[
    \rule{0em}{2ex}
    u(1,\theta_{1})
    +
    \ii
    v(1,\theta_{1})
  \right]
  \left[
    1
    -
    \ii\,
    \frac{\cos(\Delta\theta/2)}{\sin(\Delta\theta/2)}
  \right]
  \nonumber\\
  & = &
  \frac{1}{2\pi}\,
  \PV\int_{-\pi}^{\pi}d\theta_{1}\,
  \left[
    \rule{0em}{2ex}
    u(1,\theta_{1})
    +
    \ii
    v(1,\theta_{1})
  \right]
  +
  \nonumber\\
  &   &
  +
  \frac{1}{2\pi}\,
  \PV\int_{-\pi}^{\pi}d\theta_{1}\,
  \frac{\cos(\Delta\theta/2)}{\sin(\Delta\theta/2)}
  \left[
    \rule{0em}{2ex}
    v(1,\theta_{1})
    -
    \ii 
    u(1,\theta_{1})
  \right].
\end{eqnarray}

\noindent
Since both $u(1,\theta_{1})$ and $v(1,\theta_{1})$ are zero-average real
functions on the unit circle, the first two integrals in the last form of
the equation above are zero, so that separating the real and imaginary
parts within this expression we are left with

\noindent
\begin{eqnarray}
  u(1,\theta)
  +
  \ii
  v(1,\theta)
  & = &
  \frac{1}{2\pi}\,
  \PV\int_{-\pi}^{\pi}d\theta_{1}\,
  \frac{\cos(\Delta\theta/2)}{\sin(\Delta\theta/2)}\,
  v(1,\theta_{1})
  +
  \nonumber\\
  &   &
  -\,
  \frac{\ii}{2\pi}\,
  \PV\int_{-\pi}^{\pi}d\theta_{1}\,
  \frac{\cos(\Delta\theta/2)}{\sin(\Delta\theta/2)}\,
  u(1,\theta_{1}),
\end{eqnarray}

\noindent
where $\Delta\theta=\theta_{1}-\theta$. Separating the real and imaginary
parts of this equation we may now write that

\noindent
\begin{eqnarray}
  u(1,\theta)
  & = &
  \frac{1}{2\pi}\,
  \PV\int_{-\pi}^{\pi}d\theta_{1}\,
  \frac{\cos(\Delta\theta/2)}{\sin(\Delta\theta/2)}\,
  v(1,\theta_{1}),
  \nonumber\\
  v(1,\theta)
  & = &
  -\,
  \frac{1}{2\pi}\,
  \PV\int_{-\pi}^{\pi}d\theta_{1}\,
  \frac{\cos(\Delta\theta/2)}{\sin(\Delta\theta/2)}\,
  u(1,\theta_{1}),
\end{eqnarray}

\noindent
where $\Delta\theta=\theta_{1}-\theta$. Recalling that
$f(\theta)=u(1,\theta)$ and $g(\theta)=v(1,\theta)$, almost everywhere on
the unit circle, we have

\noindent
\begin{eqnarray}
  f(\theta)
  & = &
  \HH^{-1}[g(\theta)]
  \nonumber\\
  & = &
  \frac{1}{2\pi}\,
  \PV\int_{-\pi}^{\pi}d\theta_{1}\,
  \frac{\cos(\Delta\theta/2)}{\sin(\Delta\theta/2)}\,
  g(\theta_{1}),
  \nonumber\\
  g(\theta)
  & = &
  \HH[f(\theta)]
  \nonumber\\
  & = &
  -\,
  \frac{1}{2\pi}\,
  \PV\int_{-\pi}^{\pi}d\theta_{1}\,
  \frac{\cos(\Delta\theta/2)}{\sin(\Delta\theta/2)}\,
  f(\theta_{1}),
\end{eqnarray}

\noindent
two equations which are thus valid almost everywhere as well. These are
the transformations relating the pair of Fourier-conjugate functions
$f(\theta)$ and $g(\theta)$. The second expression defines the compact
Hilbert transformation of $f(\theta)$ into $g(\theta)$, and the first one
defines the inverse transformation, which recovers $f(\theta)$ from
$g(\theta)$. Note that in this notation the transform is defined with an
explicit minus sign, and that its inverse is simply minus the transform
itself, $\HH^{-1}[g(\theta)]=-\HH[g(\theta)]$. This completes the proof of
Theorem~\ref{Theo01}.

\vspace{2.6ex}

\noindent
It is interesting to observe that this transform can be interpreted as a
linear integral operator acting on the space of zero-average integrable
real functions defined on the unit circle. The integration kernel of the
integral operator depends only on the difference $\theta-\theta_{1}$, and
is given by

\begin{equation}\label{EQIntKernHc}
  K_{\HH}(\theta-\theta_{1})
  =
  -\,
  \frac{1}{2\pi}\,
  \frac
  {\cos\!\left[\rule{0em}{2ex}(\theta_{1}-\theta)/2\right]}
  {\sin\!\left[\rule{0em}{2ex}(\theta_{1}-\theta)/2\right]},
\end{equation}

\noindent
so that the action of the operator on an arbitrarily given zero-average
real integrable function $f(\theta)$ on the unit circle can be written as

\noindent
\begin{eqnarray}
  g(\theta)
  & = &
  \HH[f(\theta)]
  \nonumber\\
  & = &
  \PV\int_{-\pi}^{\pi}d\theta_{1}\,
  K_{\HH}(\theta-\theta_{1})
  f(\theta_{1}).
\end{eqnarray}

\noindent
The operator is linear, invertible, and the composition of the operator
with itself results in the operation of multiplication by $-1$. Note that,
since by hypothesis $f(\theta)$ is integrable on the unit circle, the
Cauchy principal value refers only to the explicit non-integrable
singularity of the integration kernel at the position $\theta_{1}=\theta$.

\section{Action on the Fourier Basis}\label{Sec04}

We will now determine the action of the compact Hilbert transform on the
elements of the Fourier basis of functions. The case of the constant
function, which constitutes the $k=0$ element of the basis, that is the
single member of the basis which is not a zero-average function, must be
examined in separate. We will now prove the following simple theorem.

\begin{theorem}\Colon\label{Theo02}
  Given any constant real function $f(\theta)=R$, for any real constant
  $R$, its compact Hilbert transform is zero, that is, $\HH[R]=0$.
\end{theorem}

\begin{proof}\Colon
\end{proof}

\noindent
We start from the expression in Equation~(\ref{EQCauchyPV}) for the very
simple case $w(z)=1$,

\begin{equation}
  1
  =
  \frac{1}{\pi\ii}\,
  \PV\oint_{C}dz_{1}\,
  \frac{1}{z_{1}-z},
\end{equation}

\noindent
where both $z_{1}$ and $z$ are on the circle $C$ of radius $\rho_{1}$. We
may now write all quantities in this equation in terms of $\rho_{1}$,
$\theta_{1}$ and $\theta$,

\noindent
\begin{eqnarray}
  1
  & = &
  \frac{1}{\pi\ii}\,
  \PV\int_{-\pi}^{\pi}d\theta_{1}\,
  \ii\rho_{1}\e{\iii\theta_{1}}\,
  \frac{1}{\rho_{1}\e{\iii\theta_{1}}-\rho_{1}\e{\iii\theta}}
  \nonumber\\
  & = &
  \frac{1}{\pi}\,
  \PV\int_{-\pi}^{\pi}d\theta_{1}\,
  \frac{1}{1-\e{-\iii\Delta\theta}},
\end{eqnarray}

\noindent
where $\Delta\theta=\theta_{1}-\theta$. Note that, since there are no
remaining dependencies on $\rho_{1}$, we may now take the $\rho_{1}\to
1_{(-)}$ limit of this expression, in which the principal value acquires
its usual real meaning on the unit circle. Just as in the previous
section, in order to identify separately the real and imaginary parts of
this equation, we must now rationalize the integrand. Using the result in
Equation~(\ref{EQRatiForm1}) we obtain

\noindent
\begin{eqnarray}
  1
  & = &
  \frac{1}{2\pi}\,
  \PV\int_{-\pi}^{\pi}d\theta_{1}\,
  \left[
    1
    -
    \ii\,
    \frac{\cos(\Delta\theta/2)}{\sin(\Delta\theta/2)}
  \right]
  \nonumber\\
  & = &
  \frac{1}{2\pi}\,
  \PV\int_{-\pi}^{\pi}d\theta_{1}\,
  -\,
  \frac{\ii}{2\pi}\,
  \PV\int_{-\pi}^{\pi}d\theta_{1}\,
  \frac{\cos(\Delta\theta/2)}{\sin(\Delta\theta/2)}
  \nonumber\\
  & = &
  1
  -\,
  \frac{\ii}{2\pi}\,
  \PV\int_{-\pi}^{\pi}d\theta_{1}\,
  \frac{\cos(\Delta\theta/2)}{\sin(\Delta\theta/2)}.
\end{eqnarray}

\noindent
It follows therefore that we have

\begin{equation}
  -\,
  \frac{1}{2\pi}\,
  \PV\int_{-\pi}^{\pi}d\theta_{1}\,
  \frac{\cos(\Delta\theta/2)}{\sin(\Delta\theta/2)}
  =
  0,
\end{equation}

\noindent
which is the statement that $\HH[1]=0$. Note that, since
$\Delta\theta=\theta_{1}-\theta$, which in the context of this integral
implies that $d\theta_{1}=d(\Delta\theta)$, by means of a trivial
transformation of variables this integral can also be shown to be zero by
simple parity arguments. Given the linearity of the compact Hilbert
transform, it is equally true that, for any real constant $R$, we have
that $\HH[R]=0$, so that all constant functions are mapped to the null
function. This completes the proof of Theorem~\ref{Theo02}.

\vspace{2.6ex}

\noindent
Let us now consider all the remaining elements of the Fourier basis of
functions. We will prove the following theorem.

\begin{theorem}\Colon\label{Theo03}
  Given the elements of the Fourier basis of functions, $\cos(k\theta)$
  and $\sin(k\theta)$, for $k\in\{1,2,3,\ldots,\infty\}$, the following
  relations between them hold:
\end{theorem}

\noindent
\begin{eqnarray}
  \cos(k\theta)
  & = &
  -
  \HH\!\left[\rule{0em}{2ex}\sin(k\theta)\right],
  \nonumber\\
  \sin(k\theta)
  & = &
  \HH\!\left[\rule{0em}{2ex}\cos(k\theta)\right].
\end{eqnarray}

\begin{proof}\Colon
\end{proof}

\noindent
In order to prove this theorem we start from the expression in
Equation~(\ref{EQCauchyPV}) for the case $w(z)=z^{k}$, where
$k\in\{1,2,3,\ldots,\infty\}$, that is, for a strictly positive power of
$z$, which is therefore an inner analytic function. Note that these are
all the elements of the complex Taylor basis of functions, with the
exception of the constant function. We have therefore

\begin{equation}
  z^{k}
  =
  \frac{1}{\pi\ii}\,
  \PV\oint_{C}dz_{1}\,
  \frac{z_{1}^{k}}{z_{1}-z},
\end{equation}

\noindent
where both $z_{1}$ and $z$ are on the circle $C$ of radius $\rho_{1}$. We
may now write all quantities in this equation in terms of $\rho_{1}$,
$\theta_{1}$ and $\theta$,

\noindent
\begin{eqnarray}
  \rho_{1}^{k}
  \e{\iii k\theta}
  & = &
  \frac{1}{\pi\ii}\,
  \PV\int_{-\pi}^{\pi}d\theta_{1}\,
  \ii\rho_{1}\e{\iii\theta_{1}}\,
  \frac
  {\rho_{1}^{k}\e{\iii k\theta_{1}}}
  {\rho_{1}\e{\iii\theta_{1}}-\rho_{1}\e{\iii\theta}}
  \;\;\;\Rightarrow
  \nonumber\\
  \e{\iii k\theta}
  & = &
  \frac{1}{\pi}\,
  \PV\int_{-\pi}^{\pi}d\theta_{1}\,
  \frac
  {\e{\iii k\theta_{1}}}
  {1-\e{-\iii\Delta\theta}},
\end{eqnarray}

\noindent
where $\Delta\theta=\theta_{1}-\theta$. Note that, since there are no
remaining dependencies on $\rho_{1}$, we may now take the $\rho_{1}\to
1_{(-)}$ limit of this expression, in which the principal value acquires
its usual real meaning on the unit circle. In the limit we have

\begin{equation}
  \cos(k\theta)
  +
  \ii
  \sin(k\theta)
  =
  \frac{1}{\pi}\,
  \PV\int_{-\pi}^{\pi}d\theta_{1}\,
  \frac
  {\cos(k\theta_{1})+\ii\sin(k\theta_{1})}
  {1-\e{-\iii\Delta\theta}}.
\end{equation}

\noindent
Just as in the previous cases, in order to identify separately the real
and imaginary parts of this equation, we must now rationalize the
integrand. Using the result in Equation~(\ref{EQRatiForm1}) we obtain

\noindent
\begin{eqnarray}
  \lefteqn
  {
    \cos(k\theta)
    +
    \ii
    \sin(k\theta)
  }
  &   &
  \nonumber\\
  & = &
  \frac{1}{2\pi}\,
  \PV\int_{-\pi}^{\pi}d\theta_{1}\,
  \left[
    \rule{0em}{2ex}
    \cos(k\theta_{1})
    +
    \ii
    \sin(k\theta_{1})
  \right]
  \left[
    1
    -
    \ii\,
    \frac{\cos(\Delta\theta/2)}{\sin(\Delta\theta/2)}
  \right]
  \nonumber\\
  & = &
  \frac{1}{2\pi}\,
  \PV\int_{-\pi}^{\pi}d\theta_{1}\,
  \left[
    \rule{0em}{2ex}
    \cos(k\theta_{1})
    +
    \ii
    \sin(k\theta_{1})
  \right]
  +
  \nonumber\\
  &   &
  +
  \frac{1}{2\pi}\,
  \PV\int_{-\pi}^{\pi}d\theta_{1}\,
  \frac{\cos(\Delta\theta/2)}{\sin(\Delta\theta/2)}
  \left[
    \rule{0em}{2ex}
    \sin(k\theta_{1})
    -
    \ii
    \cos(k\theta_{1})
  \right],
\end{eqnarray}

\noindent
where $\Delta\theta=\theta_{1}-\theta$ and $k\in\{1,2,3,\ldots,\infty\}$.
The first two integrals in the last form of the equation above are zero
for all $k>0$ because they are integrals of cosines and sines over integer
multiples of their periods, so that we may now separate the real and
imaginary parts of the remaining terms and thus get

\noindent
\begin{eqnarray}
  \cos(k\theta)
  +
  \ii
  \sin(k\theta)
  & = &
  \frac{1}{2\pi}\,
  \PV\int_{-\pi}^{\pi}d\theta_{1}\,
  \frac{\cos(\Delta\theta/2)}{\sin(\Delta\theta/2)}\,
  \sin(k\theta_{1})
  +
  \nonumber\\
  &   &
  -\,
  \frac{\ii}{2\pi}\,
  \PV\int_{-\pi}^{\pi}d\theta_{1}\,
  \frac{\cos(\Delta\theta/2)}{\sin(\Delta\theta/2)}\,
  \cos(k\theta_{1}),
\end{eqnarray}

\noindent
where $\Delta\theta=\theta_{1}-\theta$ and $k\in\{1,2,3,\ldots,\infty\}$.
We therefore obtain the action of the compact Hilbert transform on the
elements of the Fourier basis,

\noindent
\begin{eqnarray}
  \cos(k\theta)
  & = &
  -
  \HH\!\left[\rule{0em}{2ex}\sin(k\theta)\right]
  \nonumber\\
  & = &
  \frac{1}{2\pi}\,
  \PV\int_{-\pi}^{\pi}d\theta_{1}\,
  \frac{\cos(\Delta\theta/2)}{\sin(\Delta\theta/2)}\,
  \sin(k\theta_{1}),
  \nonumber\\
  \sin(k\theta)
  & = &
  \HH\!\left[\rule{0em}{2ex}\cos(k\theta)\right]
  \nonumber\\
  & = &
  -\,
  \frac{1}{2\pi}\,
  \PV\int_{-\pi}^{\pi}d\theta_{1}\,
  \frac{\cos(\Delta\theta/2)}{\sin(\Delta\theta/2)}\,
  \cos(k\theta_{1}),
\end{eqnarray}

\noindent
where $\Delta\theta=\theta_{1}-\theta$ and $k\in\{1,2,3,\ldots,\infty\}$.
The second equation above is the transform applied to the cosines and the
first equation is the inverse transform applied to the sines. As one can
see, the transform does indeed have the property of replacing cosines with
sines and sines with minus cosines, as expected. This completes the proof
of Theorem~\ref{Theo03}.

\vspace{2.6ex}

\noindent
One can now see that the application of the compact Hilbert transform to
the Fourier series of an arbitrarily given zero-average integrable real
function $f(\theta)$ on the unit circle will produce the Fourier series of
its Fourier-conjugate real function $g(\theta)$. Given the linearity of
the transform, if we apply it to the Fourier series of $f(\theta)$ we get

\noindent
\begin{eqnarray}
  \HH\!\left\{
    \sum_{k=1}^{\infty}
    \left[
      \rule{0em}{2ex}
      \alpha_{k}\cos(k\theta)
      +
      \beta_{k}\sin(k\theta)
    \right]
  \right\}
  & = &
  \sum_{k=1}^{\infty}
  \left\{
    \rule{0em}{2.5ex}
    \alpha_{k}
    \HH\!\left[\rule{0em}{2ex}\cos(k\theta)\right]
    +
    \beta_{k}
    \HH\!\left[\rule{0em}{2ex}\sin(k\theta)\right]
  \right\}
  \nonumber\\
  & = &
  \sum_{k=1}^{\infty}
  \left[
    \rule{0em}{2ex}
    \alpha_{k}
    \sin(k\theta)
    -
    \beta_{k}
    \cos(k\theta)
  \right],
\end{eqnarray}

\noindent
where this last one is the Fourier series of the real function
$g(\theta)$, which is the Fourier conjugate of $f(\theta)$. Hence, if
$S(\rho,\theta)$ is the complex power series given in
Equation~(\ref{EQPowrSers}), if $S^{F,f}(\theta)=\Re[S(1,\theta)]$ is the
Fourier series of $f(\theta)$ and $S^{F,g}(\theta)=\Im[S(1,\theta)]$ is
the Fourier series of $g(\theta)$, then we have that

\begin{equation}
  S^{F,g}(\theta)
  =
  \HH\!\left[\rule{0em}{2.2ex}S^{F,f}(\theta)\right].
\end{equation}

\noindent
Note that the same is true for the corresponding partial sums, as well as
for the corresponding remainders, so long as the latter exist at all. If
$S_{N}(\rho,\theta)$ is the $N^{\rm th}$ partial sum and
$R_{N}(\rho,\theta)$ is the $N^{\rm th}$ remainder of the complex power
series given in Equation~(\ref{EQPowrSers}), and if
$S_{N}^{F,f}(\theta)=\Re[S_{N}(1,\theta)]$ is the $N^{\rm th}$ partial sum
of the real Fourier series of $f(\theta)$, if
$S_{N}^{F,g}(\theta)=\Im[S_{N}(1,\theta)]$ is the $N^{\rm th}$ partial sum
of the Fourier series of $g(\theta)$, if
$R_{N}^{F,f}(\theta)=\Re[R_{N}(1,\theta)]$ is the $N^{\rm th}$ remainder
of the Fourier series of $f(\theta)$ and if
$R_{N}^{F,g}(\theta)=\Im[R_{N}(1,\theta)]$ is the $N^{\rm th}$ remainder
of the Fourier series of $g(\theta)$, then we have

\noindent
\begin{eqnarray}
  S_{N}^{F,g}(\theta)
  & = &
  \HH\!\left[\rule{0em}{2ex}S_{N}^{F,f}(\theta)\right],
  \nonumber\\\label{EQRemMap}
  R_{N}^{F,g}(\theta)
  & = &
  \HH\!\left[\rule{0em}{2ex}R_{N}^{F,f}(\theta)\right].
\end{eqnarray}

\noindent
Of course, in each one of these cases the inverse mapping holds as well,
using the inverse transform to take us from the quantities related to
$g(\theta)$ back to the corresponding quantities related to $f(\theta)$.

\section{An Infinite Collection of Identities}\label{Sec05}

In order to obtain a certain infinite collection of identities satisfied
by all zero-average integrable real functions and their Fourier-conjugate
real functions, which will be very important later, we start by examining
the action of the compact Hilbert transform on the products of arbitrarily
given integrable real functions and the elements of the Fourier basis. We
will prove the following theorem.

\begin{theorem}\Colon\label{Theo04}
  Given an arbitrary zero-average integrable real function $f(\theta)$ on
  the unit circle, and the corresponding Fourier-conjugate real function
  $g(\theta)$, these two real functions satisfy almost everywhere the
  following infinite collection of identities:
\end{theorem}

\noindent
\begin{eqnarray}\label{EQTheo04}
  f(\theta)
  & = &
  \frac{1}{2\pi}\,
  \PV\int_{-\pi}^{\pi}d\theta_{1}\,
  \frac
  {
    \sin\!\left[\rule{0em}{2ex}(k+1/2)\Delta\theta\right]
    f(\theta_{1})
    +
    \cos\!\left[\rule{0em}{2ex}(k+1/2)\Delta\theta\right]
    g(\theta_{1})
  }
  {
    \sin(\Delta\theta/2)
  },
  \nonumber\\
  g(\theta)
  & = &
  \frac{1}{2\pi}\,
  \PV\int_{-\pi}^{\pi}d\theta_{1}\,
  \frac
  {
    \sin\!\left[\rule{0em}{2ex}(k+1/2)\Delta\theta\right]
    g(\theta_{1})
    -
    \cos\!\left[\rule{0em}{2ex}(k+1/2)\Delta\theta\right]
    f(\theta_{1})
  }
  {
    \sin(\Delta\theta/2)
  },
\end{eqnarray}

\noindent
where $\Delta\theta=\theta_{1}-\theta$ and $k\in\{1,2,3,\ldots,\infty\}$.

\begin{proof}\Colon
\end{proof}

\noindent
In order to prove this theorem we start from the expression in
Equation~(\ref{EQCauchyPV}), exchanging $w(z)$ for the product
$z^{k}w(z)$, which is also an inner analytic function so long as $z^{k}$
is an arbitrary positive integer power, which it is since we assume that
$k\in\{1,2,3,\ldots,\infty\}$. We therefore have

\begin{equation}
  z^{k}
  w(z)
  =
  \frac{1}{\pi\ii}\,
  \PV\oint_{C}dz_{1}\,
  \frac{z_{1}^{k}w(z_{1})}{z_{1}-z},
\end{equation}

\noindent
where both $z_{1}$ and $z$ are on the circle $C$ of radius $\rho_{1}$
within the open unit disk. We may now write all quantities in this
equation in terms of $\rho_{1}$, $\theta_{1}$ and $\theta$,

\noindent
\begin{eqnarray}
  \rho_{1}^{k}
  \e{\iii k\theta}
  w(\rho_{1},\theta)
  & = &
  \frac{1}{\pi\ii}\,
  \PV\int_{-\pi}^{\pi}d\theta_{1}\,
  \ii\rho_{1}\e{\iii\theta_{1}}\,
  \frac
  {\rho_{1}^{k}\e{\iii k\theta_{1}}w(\rho_{1},\theta_{1})}
  {\rho_{1}\e{\iii\theta_{1}}-\rho_{1}\e{\iii\theta}}
  \;\;\;\Rightarrow
  \nonumber\\
  w(\rho_{1},\theta)
  & = &
  \frac{1}{\pi}\,
  \PV\int_{-\pi}^{\pi}d\theta_{1}\,
  \frac
  {\e{\iii k\Delta\theta}w(\rho_{1},\theta_{1})}
  {1-\e{-\iii\Delta\theta}},
\end{eqnarray}

\noindent
where $\Delta\theta=\theta_{1}-\theta$ and $k\in\{1,2,3,\ldots,\infty\}$.
Note that, since by construction the real and imaginary parts
$u(\rho_{1},\theta)$ and $v(\rho_{1},\theta)$ of $w(\rho_{1},\theta)$ for
$\rho_{1}=1$ are integrable real functions on the unit circle, and since
there are no other dependencies on $\rho_{1}$ in this equation, we may now
take the $\rho_{1}\to 1_{(-)}$ limit of this expression, in which the
principal value acquires its usual real meaning on the unit circle, thus
obtaining

\begin{equation}
  u(1,\theta)
  +
  \ii
  v(1,\theta)
  =
  \frac{1}{\pi}\,
  \PV\int_{-\pi}^{\pi}d\theta_{1}\,
  \frac
  {
    \e{\iii k\Delta\theta}
    \left[
      \rule{0em}{2ex}
      u(1,\theta_{1})
      +
      \ii
      v(1,\theta_{1})
    \right]
  }
  {
    1-\e{-\iii\Delta\theta}
  },
\end{equation}

\noindent
where $\Delta\theta=\theta_{1}-\theta$ and $k\in\{1,2,3,\ldots,\infty\}$.
Once more, in order to identify separately the real and imaginary parts of
this equation, we must now rationalize the integrand. Using this time the
form shown in Equation~(\ref{EQRatiForm2}) for the factor to be
rationalized, we get

\noindent
\begin{eqnarray}\label{EQInfColIdenI}
  \lefteqn
  {
    u(1,\theta)
    +
    \ii
    v(1,\theta)
  }
  &   &
  \nonumber\\
  & = &
  \frac{1}{2\pi}\,
  \PV\int_{-\pi}^{\pi}d\theta_{1}\,
  \e{\iii k\Delta\theta}\,
  \left[
    \rule{0em}{2ex}
    u(1,\theta_{1})
    +
    \ii
    v(1,\theta_{1})
  \right]
  (-\ii)\,
  \frac
  {\e{\iii\Delta\theta/2}}
  {\sin(\Delta\theta/2)}
  \nonumber\\
  & = &
  \frac{1}{2\pi}\,
  \PV\int_{-\pi}^{\pi}d\theta_{1}\,
  \left[
    \rule{0em}{2ex}
    v(1,\theta_{1})
    -
    \ii
    u(1,\theta_{1})
  \right]
  \frac
  {\e{\iii(k+1/2)\Delta\theta}}
  {\sin(\Delta\theta/2)}
  \nonumber\\
  & = &
  \frac{1}{2\pi}\,
  \PV\int_{-\pi}^{\pi}d\theta_{1}\,
  \frac
  {
    \left[
      \rule{0em}{2ex}
      v(1,\theta_{1})
      -
      \ii
      u(1,\theta_{1})
    \right]
    \left[
      \rule{0em}{2ex}
      \cos(k_{1}\Delta\theta)
      +
      \ii
      \sin(k_{1}\Delta\theta)
    \right]
  }
  {
    \sin(\Delta\theta/2)
  },
\end{eqnarray}

\noindent
where $\Delta\theta=\theta_{1}-\theta$ and $k_{1}=k+1/2$, with
$k\in\{1,2,3,\ldots,\infty\}$. Expanding the numerator in the integrand of
this integral we have

\noindent
\begin{eqnarray}
  \lefteqn
  {
    \left[
      \rule{0em}{2ex}
      v(1,\theta_{1})
      -
      \ii
      u(1,\theta_{1})
    \right]
    \left[
      \rule{0em}{2ex}
      \cos(k_{1}\Delta\theta)
      +
      \ii
      \sin(k_{1}\Delta\theta)
    \right]
  }
  &   &
  \nonumber\\
  & = &
  \left[
    \rule{0em}{2ex}
    \sin(k_{1}\Delta\theta)
    u(1,\theta_{1})
    +
    \cos(k_{1}\Delta\theta)
    v(1,\theta_{1})
  \right]
  +
  \nonumber\\
  &   &
  +
  \ii
  \left[
    \rule{0em}{2ex}
    \sin(k_{1}\Delta\theta)
    v(1,\theta_{1})
    -
    \cos(k_{1}\Delta\theta)
    u(1,\theta_{1})
  \right],
\end{eqnarray}

\noindent
and therefore we are left with

\noindent
\begin{eqnarray}
  \lefteqn
  {
    u(1,\theta)
    +
    \ii
    v(1,\theta)
  }
  &   &
  \nonumber\\
  & = &
  \frac{1}{2\pi}\,
  \PV\int_{-\pi}^{\pi}d\theta_{1}\,
  \frac
  {
    \sin(k_{1}\Delta\theta)
    u(1,\theta_{1})
    +
    \cos(k_{1}\Delta\theta)
    v(1,\theta_{1})
  }
  {\sin(\Delta\theta/2)}
  \nonumber\\
  &   &
  +
  \ii\,
  \frac{1}{2\pi}\,
  \PV\int_{-\pi}^{\pi}d\theta_{1}\,
  \frac
  {
    \sin(k_{1}\Delta\theta)
    v(1,\theta_{1})
    -
    \cos(k_{1}\Delta\theta)
    u(1,\theta_{1})
  }
  {\sin(\Delta\theta/2)},
\end{eqnarray}

\noindent
where $\Delta\theta=\theta_{1}-\theta$ and $k_{1}=k+1/2$, with
$k\in\{1,2,3,\ldots,\infty\}$. Separating the real and imaginary parts we
therefore obtain an infinite collection of identities in the form

\noindent
\begin{eqnarray}\label{EQInfColIdenII}
  \lefteqn
  {
    u(1,\theta)
  }
  &   &
  \nonumber\\
  & = &
  \frac{1}{2\pi}\,
  \PV\int_{-\pi}^{\pi}d\theta_{1}\,
  \frac
  {
    \sin\!\left[\rule{0em}{2ex}(k+1/2)\Delta\theta\right]
    u(1,\theta_{1})
    +
    \cos\!\left[\rule{0em}{2ex}(k+1/2)\Delta\theta\right]
    v(1,\theta_{1})
  }
  {
    \sin(\Delta\theta/2)
  },
  \nonumber\\
  \lefteqn
  {
    v(1,\theta)
  }
  &   &
  \nonumber\\
  & = &
  \frac{1}{2\pi}\,
  \PV\int_{-\pi}^{\pi}d\theta_{1}\,
  \frac
  {
    \sin\!\left[\rule{0em}{2ex}(k+1/2)\Delta\theta\right]
    v(1,\theta_{1})
    -
    \cos\!\left[\rule{0em}{2ex}(k+1/2)\Delta\theta\right]
    u(1,\theta_{1})
  }
  {
    \sin(\Delta\theta/2)
  },
\end{eqnarray}

\noindent
where $\Delta\theta=\theta_{1}-\theta$ and $k\in\{1,2,3,\ldots,\infty\}$.
Recalling now that $f(\theta)=u(1,\theta)$, and also that
$g(\theta)=v(1,\theta)$, almost everywhere over the unit circle, one
obtains the results in Equation~(\ref{EQTheo04}), and therefore this
completes the proof of Theorem~\ref{Theo04}.

\vspace{2.6ex}

\noindent
Note that since this is an infinite collection of integral identities,
satisfied by $f(\theta)$ and $g(\theta)$ for all strictly positive $k$, it
follows that the right-hand sides of the equations above do not, in fact,
depend on $k$. If one recognizes in the first term of each one of these
two equations the well-known result for the $k^{\rm th}$ partial sums of
the corresponding Fourier series in terms of Dirichlet
integrals~\cite{FSchurchill}, then it follows that the other terms must be
the corresponding remainders. This provides us with some level of
understanding of the nature of this infinite set of identities. In the
next section we will prove that one does obtain in fact the partial sums
and remainders of the corresponding Fourier series directly from our
complex-analytic structure.

\section{Remainders of Fourier Series}\label{Sec06}

We will now derive certain expressions for the partial sums and for the
corresponding remainders of the Fourier series. In order to do this, let
$f(\theta)$ be a zero-average integrable real function defined on
$[-\pi,\pi]$ and let the real numbers $\alpha_{0}=0$, $\alpha_{k}$ and
$\beta_{k}$, for $k\in\{1,2,3,\ldots,\infty\}$, be its Fourier
coefficients. We then define the complex coefficients $c_{0}=0$ and
$c_{k}$ shown in Equation~(\ref{EQTaylCoef}), and thus construct the
corresponding proper inner analytic function $w(z)$ within the open unit
disk, using the power series $S(z)$ given in Equation~(\ref{EQPowrSers}),
which, as was shown in~\cite{CAoRFI}, always converges for
$|z|<1$. Considering that $c_{0}\equiv 0$, the partial sums of the first
$N$ terms of this series are given by

\begin{equation}\label{EQPowSerPar}
  S_{N}(z)
  =
  \sum_{k=0}^{N}
  c_{k}z^{k},
\end{equation}

\noindent
where $N\in\{1,2,3,\ldots,\infty\}$, a complex sequence for each value of
$z$ which, for $|z|<1$, we already know to converge to $w(z)$ in the
$N\to\infty$ limit. Note however that, since $S_{N}(z)$ is a polynomial of
order $N$ and therefore an analytic function over the whole complex plane,
this expression itself can be consistently considered for all finite $N$
and all $z$, and in particular for $z$ on the unit circle, where $|z|=1$.
One can also define the corresponding remainders of the complex power
series, in the usual way, as

\begin{equation}\label{EQPowSerRem}
  R_{N}(z)
  =
  w(z)
  -
  S_{N}(z).
\end{equation}

\noindent
We will now prove the following theorem.

\begin{theorem}\Colon\label{Theo05}
  Given an arbitrary zero-average integrable real function $f(\theta)$ on
  the unit circle and the corresponding Fourier-conjugate real function
  $g(\theta)$, if $S_{N}^{F}(\theta)=\Re[S_{N}(1,\theta)]$ is the $N^{\rm
    th}$ partial sum of the Fourier series of $f(\theta)$, and if
  $R_{N}^{F}(\theta)=\Re[R_{N}(1,\theta)]$ is the corresponding remainder
  of that Fourier series, then we have that this partial sum and this
  remainder are given by the following integrals:
\end{theorem}

\noindent
\begin{eqnarray}
  S_{N}^{F}(\theta)
  & = &
  \frac{1}{2\pi}\,
  \int_{-\pi}^{\pi}d\theta_{1}\,
  \frac
  {
    \sin\!\left[\rule{0em}{2ex}(N+1/2)\Delta\theta\right]
  }
  {
    \sin(\Delta\theta/2)
  }\,
  f(\theta_{1}),
  \nonumber\\
  R_{N}^{F}(\theta)
  & = &
  \frac{1}{2\pi}\,
  \PV\int_{-\pi}^{\pi}d\theta_{1}\,
  \frac
  {
    \cos\!\left[\rule{0em}{2ex}(N+1/2)\Delta\theta\right]
  }
  {
    \sin(\Delta\theta/2)
  }\,
  g(\theta_{1}),
\end{eqnarray}

\noindent
where $\Delta\theta=\theta_{1}-\theta$ and $N\in\{1,2,3,\ldots,\infty\}$.

\vspace{2.6ex}

\noindent
Note that the integral in the expression of the partial sum is the known
Dirichlet integral, while the one in the expression of the reminder is
similar but not identical to it.

\begin{proof}\Colon
\end{proof}

\noindent
In order to prove this theorem, let us consider the complex partial sums
$S_{N}(z)$ as given in Equation~(\ref{EQPowSerPar}). In addition to this,
the complex coefficients $c_{k}$ may be written as integrals involving
$w(z)$, with the use of the Cauchy integral formulas,

\begin{equation}
  c_{k}
  =
  \frac{1}{2\pi\ii}
  \oint_{C}dz\,
  \frac{w(z)}{z^{k+1}},
\end{equation}

\noindent
for $k\in\{0,1,2,3,\ldots,\infty\}$, where $C$ can be taken as a circle
centered at the origin, with radius $\rho\leq 1$. The reason why we may
include the case $\rho=1$ here is that, as was shown in~\cite{CAoRFI}, as
a function of $\rho$ the expression above for $c_{k}$ is not only constant
within the open unit disk, but also continuous from within at the unit
circle. In this way the coefficients $c_{k}$ may be written back in terms
of the inner analytic function $w(z)$. If we substitute this expression
for $c_{k}$ back in the partial sums of the complex power series shown in
Equation~(\ref{EQPowSerPar}) we get

\noindent
\begin{eqnarray}
  S_{N}(z)
  & = &
  \sum_{k=0}^{N}
  z^{k}\,
  \frac{1}{2\pi\ii}
  \oint_{C}dz_{1}\,
  \frac{w(z_{1})}{z_{1}^{k+1}}
  \nonumber\\
  & = &
  \frac{1}{2\pi\ii}
  \oint_{C}dz_{1}\,
  \frac{w(z_{1})}{z_{1}}
  \sum_{k=0}^{N}
  \left(
    \frac{z}{z_{1}}
  \right)^{k},
\end{eqnarray}

\noindent
where $z$ can have any value, but where we must have $|z_{1}|\leq 1$. The
sum is now a finite geometric progression, so that we have its value in
closed form,

\noindent
\begin{eqnarray}
  S_{N}(z)
  & = &
  \frac{1}{2\pi\ii}
  \oint_{C}dz_{1}\,
  \frac{w(z_{1})}{z_{1}}\,
  \frac
  {1-(z/z_{1})^{N+1}}
  {1-(z/z_{1})}
  \nonumber\\
  & = &
  \frac{1}{2\pi\ii}
  \oint_{C}dz_{1}\,
  \frac{w(z_{1})}{z_{1}-z}
  -
  \frac{z^{N+1}}{2\pi\ii}
  \oint_{C}dz_{1}\,
  \frac{w(z_{1})}{z_{1}^{N+1}(z_{1}-z)}.
\end{eqnarray}

\noindent
There are now two relevant cases to be considered here, the case in which
$|z|<|z_{1}|$ and the case in which $|z|>|z_{1}|$. In the first case,
since the explicit simple pole of the integrand at the position $z_{1}=z$
lies within the integration contour, we have in the first term the Cauchy
integral formula for $w(z)$, and therefore we get

\begin{equation}\label{EQPartSums}
  S_{N}(z)
  =
  w(z)
  -
  \frac{z^{N+1}}{2\pi\ii}
  \oint_{C}dz_{1}\,
  \frac{w(z_{1})}{z_{1}^{N+1}(z_{1}-z)}.
\end{equation}

\noindent
This is the equation that allows us to write an explicit expression for
the remainder of the complex power series within the open unit disk, thus
making it easier to discuss its convergence there. In the other case, in
which $|z|>|z_{1}|$, the explicit simple pole of the integrand at the
position $z$ lies outside of the integration contour, and therefore by the
Cauchy-Goursat theorem we just have zero in the first term, so that we get

\begin{equation}
  S_{N}(z)
  =
  -\,
  \frac{z^{N+1}}{2\pi\ii}
  \oint_{C}dz_{1}\,
  \frac{w(z_{1})}{z_{1}^{N+1}(z_{1}-z)}.
\end{equation}

\noindent
This provides us, therefore, with an explicit expression for the partial
sums, but not for the remainder. The only other possible case is that in
which $|z|=|z_{1}|$, in which both $z_{1}$ and $z$ are over the circle $C$
of radius $\rho_{1}$, and therefore so is the explicit simple pole of the
integrand at the position $z_{1}=z$. In this case, just as we did before
in Section~\ref{Sec03}, we may slightly deform the integration contour $C$
in order to have it pass on one side or the other of the simple pole of
the integrand at $z_{1}=z$. If we use a deformed contour $C_{\ominus}$
that {\em excludes} the pole from its interior, then we have, instead of
Equation~(\ref{EQPartSums}),

\begin{equation}\label{EQExclPole2}
  S_{N}(z)
  =
  0
  -
  \frac{z^{N+1}}{2\pi\ii}
  \oint_{C_{\ominus}}dz_{1}\,
  \frac{w(z_{1})}{z_{1}^{N+1}(z_{1}-z)},
\end{equation}

\noindent
while if we use a deformed contour $C_{\oplus}$ that {\em includes} the
pole in its interior, then we have, just as in
Equation~(\ref{EQPartSums}),

\begin{equation}\label{EQInclPole2}
  S_{N}(z)
  =
  w(z)
  -
  \frac{z^{N+1}}{2\pi\ii}
  \oint_{C_{\oplus}}dz_{1}\,
  \frac{w(z_{1})}{z_{1}^{N+1}(z_{1}-z)}.
\end{equation}

\noindent
Once more, since by the Sokhotskii-Plemelj theorem~\cite{sokhplem} the
Cauchy principal value of the integral over $C$ is the arithmetic average
of these two integrals, taking the average of
Equations~(\ref{EQExclPole2}) and~(\ref{EQInclPole2}) we obtain the
expression

\begin{equation}
  S_{N}(z)
  =
  \frac{w(z)}{2}\,
  -
  \frac{z^{N+1}}{2\pi\ii}\,
  \PV\oint_{C}dz_{1}\,
  \frac{w(z_{1})}{z_{1}^{N+1}(z_{1}-z)},
\end{equation}

\noindent
where both $z_{1}$ and $z$ are now on the circle $C$ of radius $\rho_{1}$.
Since we have that the corresponding remainder of the series is defined as
given in Equation~(\ref{EQPowSerRem}), we get a corresponding expression
for the remainder, in terms of the same integral,

\begin{equation}
  R_{N}(z)
  =
  \frac{w(z)}{2}\,
  +
  \frac{z^{N+1}}{2\pi\ii}\,
  \PV\oint_{C}dz_{1}\,
  \frac{w(z_{1})}{z_{1}^{N+1}(z_{1}-z)},
\end{equation}

\noindent
where both $z_{1}$ and $z$ are on the circle $C$ of radius $\rho_{1}$. We
have therefore the pair of equations

\noindent
\begin{eqnarray}
  S_{N}(z)
  & = &
  \frac{w(z)}{2}\,
  -
  I_{N}(z),
  \nonumber\\
  R_{N}(z)
  & = &
  \frac{w(z)}{2}\,
  +
  I_{N}(z),
\end{eqnarray}

\noindent
and we must now write the integral $I_{N}(z)$ explicitly in terms of
$\rho_{1}$, $\theta_{1}$ and $\theta$,

\noindent
\begin{eqnarray}
  I_{N}(z)
  & = &
  \frac{z^{N+1}}{2\pi\ii}\,
  \PV\oint_{C}dz_{1}\,
  \frac{w(z_{1})}{z_{1}^{N+1}(z_{1}-z)}
  \nonumber\\
  & = &
  \frac{\rho_{1}^{N+1}\e{\iii(N+1)\theta}}{2\pi\ii}\,
  \PV\int_{-\pi}^{\pi}d\theta_{1}\,
  \ii\rho_{1}\e{\iii\theta_{1}}\,
  \frac
  {u(\rho_{1},\theta_{1})+\ii v(\rho_{1},\theta_{1})}
  {
    \rho_{1}^{N+1}
    \e{\iii(N+1)\theta_{1}}
    \left(
      \rho_{1}
      \e{\iii\theta_{1}}
      -
      \rho_{1}
      \e{\iii\theta}
    \right)
  }
  \nonumber\\
  & = &
  \frac{1}{2\pi}\,
  \PV\int_{-\pi}^{\pi}d\theta_{1}\,
  \e{-\iii(N+1)\Delta\theta}\,
  \frac
  {u(\rho_{1},\theta_{1})+\ii v(\rho_{1},\theta_{1})}
  {1-\e{-\iii\Delta\theta}},
\end{eqnarray}

\noindent
where $\Delta\theta=\theta_{1}-\theta$ and $N\in\{1,2,3,\ldots,\infty\}$.
Once again we must rationalize the integrand, and using once more the
result shown in Equation~(\ref{EQRatiForm2}) we get

\noindent
\begin{eqnarray}
  \lefteqn
  {
    I_{N}(\rho_{1},\theta)
  }
  &   &
  \nonumber\\
  & = &
  \frac{1}{4\pi}\,
  \PV\int_{-\pi}^{\pi}d\theta_{1}\,
  \e{-\iii(N+1)\Delta\theta}
  \left[
    \rule{0em}{2ex}
    u(\rho_{1},\theta_{1})
    +
    \ii
    v(\rho_{1},\theta_{1})
  \right]
  (-\ii)
  \frac
  {\e{\iii\Delta\theta/2}}
  {\sin(\Delta\theta/2)}
  \nonumber\\
  & = &
  \frac{1}{4\pi}\,
  \PV\int_{-\pi}^{\pi}d\theta_{1}\,
  \left[
    \rule{0em}{2ex}
    v(\rho_{1},\theta_{1})
    -
    \ii
    u(\rho_{1},\theta_{1})
  \right]
  \frac
  {\e{-\iii(N+1/2)\Delta\theta}}
  {\sin(\Delta\theta/2)}
  \nonumber\\
  & = &
  \frac{1}{4\pi}\,
  \PV\int_{-\pi}^{\pi}d\theta_{1}\,
  \frac
  {
    \left[
      \rule{0em}{2ex}
      v(\rho_{1},\theta_{1})
      -
      \ii
      u(\rho_{1},\theta_{1})
    \right]
    \left[
      \rule{0em}{2ex}
      \cos(N_{1}\Delta\theta)
      -
      \ii
      \sin(N_{1}\Delta\theta)
    \right]
  }
  {\sin(\Delta\theta/2)},
\end{eqnarray}

\noindent
where $\Delta\theta=\theta_{1}-\theta$ and $N_{1}=N+1/2$, with
$N\in\{1,2,3,\ldots,\infty\}$. Expanding the numerator in the integrand of
this integral we have

\noindent
\begin{eqnarray}
  \lefteqn
  {
    \left[
      \rule{0em}{2ex}
      v(\rho_{1},\theta_{1})
      -
      \ii
      u(\rho_{1},\theta_{1})
    \right]
    \left[
      \rule{0em}{2ex}
      \cos(N_{1}\Delta\theta)
      -
      \ii
      \sin(N_{1}\Delta\theta)
    \right]
  }
  &   &
  \nonumber\\
  & = &
  -
  \left[
    \rule{0em}{2ex}
    \sin(N_{1}\Delta\theta)
    u(\rho_{1},\theta_{1})
    -
    \cos(N_{1}\Delta\theta)
    v(\rho_{1},\theta_{1})
  \right]
  +
  \nonumber\\
  &   &
  -
  \ii
  \left[
    \rule{0em}{2ex}
    \sin(N_{1}\Delta\theta)
    v(\rho_{1},\theta_{1})
    +
    \cos(N_{1}\Delta\theta)
    u(\rho_{1},\theta_{1})
  \right],
\end{eqnarray}

\noindent
and therefore we are left with the following expression for our integral,

\noindent
\begin{eqnarray}
  \lefteqn
  {
    I_{N}(\rho_{1},\theta)
  }
  &   &
  \nonumber\\
  & = &
  -\,
  \frac{1}{4\pi}\,
  \PV\int_{-\pi}^{\pi}d\theta_{1}\,
  \frac
  {
    \sin(N_{1}\Delta\theta)
    u(\rho_{1},\theta_{1})
    -
    \cos(N_{1}\Delta\theta)
    v(\rho_{1},\theta_{1})
  }
  {\sin(\Delta\theta/2)}
  +
  \nonumber\\
  &   &
  -\,
  \frac{\ii}{4\pi}\,
  \PV\int_{-\pi}^{\pi}d\theta_{1}\,
  \frac
  {
    \sin(N_{1}\Delta\theta)
    v(\rho_{1},\theta_{1})
    +
    \cos(N_{1}\Delta\theta)
    u(\rho_{1},\theta_{1})
  }
  {\sin(\Delta\theta/2)},
\end{eqnarray}

\noindent
where $\Delta\theta=\theta_{1}-\theta$ and $N_{1}=N+1/2$, with
$N\in\{1,2,3,\ldots,\infty\}$. Once again we note that, since by
construction the real and imaginary parts $u(\rho_{1},\theta)$ and
$v(\rho_{1},\theta)$ of $w(\rho_{1},\theta)$ for $\rho_{1}=1$ are
integrable real functions on the unit circle, and since there are no other
dependencies on $\rho_{1}$ in this equation, we may now take the
$\rho_{1}\to 1_{(-)}$ limit of this expression, in which the principal
value acquires its usual real meaning on the unit circle, thus obtaining

\noindent
\begin{eqnarray}
  \lefteqn
  {
    I_{N}(1,\theta)
  }
  &   &
  \nonumber\\
  & = &
  -\,
  \frac{1}{4\pi}\,
  \PV\int_{-\pi}^{\pi}d\theta_{1}\,
  \frac
  {
    \sin\!\left[\rule{0em}{2ex}(N+1/2)\Delta\theta\right]
    u(1,\theta_{1})
    -
    \cos\!\left[\rule{0em}{2ex}(N+1/2)\Delta\theta\right]
    v(1,\theta_{1})
  }
  {\sin(\Delta\theta/2)}
  +
  \nonumber\\
  &   &
  -\,
  \frac{\ii}{4\pi}\,
  \PV\int_{-\pi}^{\pi}d\theta_{1}\,
  \frac
  {
    \sin\!\left[\rule{0em}{2ex}(N+1/2)\Delta\theta\right]
    v(1,\theta_{1})
    +
    \cos\!\left[\rule{0em}{2ex}(N+1/2)\Delta\theta\right]
    u(1,\theta_{1})
  }
  {\sin(\Delta\theta/2)},
\end{eqnarray}

\noindent
where $\Delta\theta=\theta_{1}-\theta$ and $N\in\{1,2,3,\ldots,\infty\}$,
in terms of which we now have the pair of equations at the unit circle

\noindent
\begin{eqnarray}
  S_{N}(1,\theta)
  & = &
  \frac{w(1,\theta)}{2}\,
  -
  I_{N}(1,\theta),
  \nonumber\\
  R_{N}(1,\theta)
  & = &
  \frac{w(1,\theta)}{2}\,
  +
  I_{N}(1,\theta).
\end{eqnarray}

\noindent
Using now the infinite collection of identities in
Equation~(\ref{EQInfColIdenII}) for the case $k=N$, which allow us to
write $w(1,\theta)/2$ in terms of integrals similar to those in
$I_{N}(1,\theta)$,

\noindent
\begin{eqnarray}
  \lefteqn
  {
    \frac{w(1,\theta)}{2}
  }
  &   &
  \nonumber\\
  & = &
  \frac{1}{4\pi}\,
  \PV\int_{-\pi}^{\pi}d\theta_{1}\,
  \frac
  {
    \sin\!\left[\rule{0em}{2ex}(N+1/2)\Delta\theta\right]
    u(1,\theta_{1})
    +
    \cos\!\left[\rule{0em}{2ex}(N+1/2)\Delta\theta\right]
    v(1,\theta_{1})
  }
  {
    \sin(\Delta\theta/2)
  }
  +
  \nonumber\\
  &   &
  +
  \frac{\ii}{4\pi}\,
  \PV\int_{-\pi}^{\pi}d\theta_{1}\,
  \frac
  {
    \sin\!\left[\rule{0em}{2ex}(N+1/2)\Delta\theta\right]
    v(1,\theta_{1})
    -
    \cos\!\left[\rule{0em}{2ex}(N+1/2)\Delta\theta\right]
    u(1,\theta_{1})
  }
  {
    \sin(\Delta\theta/2)
  },
\end{eqnarray}

\noindent
where $\Delta\theta=\theta_{1}-\theta$ and $N\in\{1,2,3,\ldots,\infty\}$,
we may write for the complex partial sums

\noindent
\begin{eqnarray}
  \lefteqn
  {
    S_{N}(1,\theta)
  }
  &   &
  \nonumber\\
  & = &
  \frac{w(1,\theta)}{2}\,
  -
  I_{N}(1,\theta)
  \nonumber\\
  & = &
  \frac{1}{4\pi}\,
  \PV\int_{-\pi}^{\pi}d\theta_{1}\,
  \frac
  {
    \sin\!\left[\rule{0em}{2ex}(N+1/2)\Delta\theta\right]
    u(1,\theta_{1})
    +
    \cos\!\left[\rule{0em}{2ex}(N+1/2)\Delta\theta\right]
    v(1,\theta_{1})
  }
  {
    \sin(\Delta\theta/2)
  }
  +
  \nonumber\\
  &   &
  +
  \frac{\ii}{4\pi}\,
  \PV\int_{-\pi}^{\pi}d\theta_{1}\,
  \frac
  {
    \sin\!\left[\rule{0em}{2ex}(N+1/2)\Delta\theta\right]
    v(1,\theta_{1})
    -
    \cos\!\left[\rule{0em}{2ex}(N+1/2)\Delta\theta\right]
    u(1,\theta_{1})
  }
  {
    \sin(\Delta\theta/2)
  }
  +
  \nonumber\\
  &   &
  +
  \frac{1}{4\pi}\,
  \PV\int_{-\pi}^{\pi}d\theta_{1}\,
  \frac
  {
    \sin\!\left[\rule{0em}{2ex}(N+1/2)\Delta\theta\right]
    u(1,\theta_{1})
    -
    \cos\!\left[\rule{0em}{2ex}(N+1/2)\Delta\theta\right]
    v(1,\theta_{1})
  }
  {\sin(\Delta\theta/2)}
  +
  \nonumber\\
  &   &
  +
  \frac{\ii}{4\pi}\,
  \PV\int_{-\pi}^{\pi}d\theta_{1}\,
  \frac
  {
    \sin\!\left[\rule{0em}{2ex}(N+1/2)\Delta\theta\right]
    v(1,\theta_{1})
    +
    \cos\!\left[\rule{0em}{2ex}(N+1/2)\Delta\theta\right]
    u(1,\theta_{1})
  }
  {\sin(\Delta\theta/2)}.
\end{eqnarray}

\noindent
As one can see in this equation, all the terms involving
$\cos\!\left[\rule{0em}{2ex}(N+1/2)\Delta\theta\right]$ cancel off, and
therefore we are left with

\noindent
\begin{eqnarray}
  S_{N}(1,\theta)
  & = &
  \frac{1}{2\pi}\,
  \PV\int_{-\pi}^{\pi}d\theta_{1}\,
  \frac
  {
    \sin\!\left[\rule{0em}{2ex}(N+1/2)\Delta\theta\right]
    u(1,\theta_{1})
  }
  {
    \sin(\Delta\theta/2)
  }
  +
  \nonumber\\
  &   &
  +
  \frac{\ii}{2\pi}\,
  \PV\int_{-\pi}^{\pi}d\theta_{1}\,
  \frac
  {
    \sin\!\left[\rule{0em}{2ex}(N+1/2)\Delta\theta\right]
    v(1,\theta_{1})
  }
  {
    \sin(\Delta\theta/2)
  },
\end{eqnarray}

\noindent
where $\Delta\theta=\theta_{1}-\theta$ and $N\in\{1,2,3,\ldots,\infty\}$.
We now observe that, since by hypothesis $f(\theta)$ is integrable on the
unit circle, the Cauchy principal value refers only to the possible
explicit non-integrable singularity of the integrands, due to the zero of
the denominators at $\theta_{1}=\theta$. However, since the numerators of
the integrands are also zero at that point, the integrands are not really
divergent at all at that point, so that from this point on we may drop the
principal value. We have therefore our final results for the real partial
sums, for both $u(1,\theta)$ and $v(1,\theta)$,

\noindent
\begin{eqnarray}
  S_{N}^{F,u}(\theta)
  & = &
  \Re[S_{N}(1,\theta)]
  \nonumber\\
  & = &
  \frac{1}{2\pi}\,
  \int_{-\pi}^{\pi}d\theta_{1}\,
  \frac
  {
    \sin\!\left[\rule{0em}{2ex}(N+1/2)\Delta\theta\right]
  }
  {
    \sin(\Delta\theta/2)
  }\,
  u(1,\theta_{1}),
  \nonumber\\
  S_{N}^{F,v}(\theta)
  & = &
  \Im[S_{N}(1,\theta)]
  \nonumber\\
  & = &
  \frac{1}{2\pi}\,
  \int_{-\pi}^{\pi}d\theta_{1}\,
  \frac
  {
    \sin\!\left[\rule{0em}{2ex}(N+1/2)\Delta\theta\right]
  }
  {
    \sin(\Delta\theta/2)
  }\,
  v(1,\theta_{1}),
\end{eqnarray}

\noindent
where $\Delta\theta=\theta_{1}-\theta$ and $N\in\{1,2,3,\ldots,\infty\}$.
These are the well-known results for the partial sums, in terms of
Dirichlet integrals~\cite{FSchurchill}. Note that the two equations above
have exactly the same form, which is to be expected, since the result
holds for all zero-average integrable real functions, including of course
both $u(1,\theta_{1})$ and $v(1,\theta_{1})$. Therefore, given an
arbitrary zero-average integrable real function $f(\theta)$ on the unit
circle, we have that the partial sums of its Fourier series are given by

\begin{equation}
  S_{N}^{F}(\theta)
  =
  \frac{1}{2\pi}\,
  \int_{-\pi}^{\pi}d\theta_{1}\,
  \frac
  {\sin\!\left[\rule{0em}{2ex}(N+1/2)(\theta_{1}-\theta)\right]}
  {\sin\!\left[\rule{0em}{2ex}(\theta_{1}-\theta)/2\right]}\,
  f(\theta_{1}),
\end{equation}

\noindent
where $N\in\{1,2,3,\ldots,\infty\}$. Note that, although this result is
already very well known, we have showed here that it does follow from our
complex-analytic structure. This completes the proof of the first part of
Theorem~\ref{Theo05}.

\vspace{2.6ex}

\noindent
Once more, it is interesting to observe that this relation can be
interpreted as a linear integral operator acting on the space of
zero-average integrable real functions defined on the unit circle, this
time resulting in the $N^{\rm th}$ partial sum of the Fourier series of a
zero-average integrable real function, a partial sum which is itself a
zero-average integrable real function. The integration kernel of this
integral operator $\DS[N,f(\theta)]$ depends only on $N$ and on the
difference $\theta-\theta_{1}$, and is given by

\begin{equation}
  K_{\DS}(N,\theta-\theta_{1})
  =
  \frac{1}{2\pi}\,
  \frac
  {\sin\!\left[\rule{0em}{2ex}(N+1/2)(\theta_{1}-\theta)\right]}
  {\sin\!\left[\rule{0em}{2ex}(\theta_{1}-\theta)/2\right]},
\end{equation}

\noindent
where $N\in\{1,2,3,\ldots,\infty\}$, so that the action of the operator on
$f(\theta)$ can be written as

\begin{equation}
  \DS[N,f(\theta)]
  =
  \int_{-\pi}^{\pi}d\theta_{1}\,
  K_{\DS}(N,\theta-\theta_{1})
  f(\theta_{1}).
\end{equation}

\noindent
Considering that its kernel is given by a Dirichlet integral, one might
call this the {\em Dirichlet operator}, so that the $N^{\rm th}$ partial
sum of the Fourier series of $f(\theta)$ is given by the action of this
operator on the zero-average integrable real function $f(\theta)$,

\begin{equation}
  S_{N}^{F}(\theta)
  =
  \DS[N,f(\theta)],
\end{equation}

\noindent
where $N\in\{1,2,3,\ldots,\infty\}$. Note that $\DS[N,f(\theta)]$
constitutes in fact a whole collection of linear integral operators acting
on the space of zero-average integrable real functions.

\begin{proof}\Colon
\end{proof}

\noindent
Using once more the very same elements that were used above for the
complex partial sums, we may also write corresponding results for the
complex remainders,

\noindent
\begin{eqnarray}
  \lefteqn
  {
    R_{N}(1,\theta)
  }
  &   &
  \nonumber\\
  & = &
  \frac{w(1,\theta)}{2}\,
  +
  I_{N}(1,\theta)
  \nonumber\\
  & = &
  \frac{1}{4\pi}\,
  \PV\int_{-\pi}^{\pi}d\theta_{1}\,
  \frac
  {
    \sin\!\left[\rule{0em}{2ex}(N+1/2)\Delta\theta\right]
    u(1,\theta_{1})
    +
    \cos\!\left[\rule{0em}{2ex}(N+1/2)\Delta\theta\right]
    v(1,\theta_{1})
  }
  {
    \sin(\Delta\theta/2)
  }
  +
  \nonumber\\
  &   &
  +
  \frac{\ii}{4\pi}\,
  \PV\int_{-\pi}^{\pi}d\theta_{1}\,
  \frac
  {
    \sin\!\left[\rule{0em}{2ex}(N+1/2)\Delta\theta\right]
    v(1,\theta_{1})
    -
    \cos\!\left[\rule{0em}{2ex}(N+1/2)\Delta\theta\right]
    u(1,\theta_{1})
  }
  {
    \sin(\Delta\theta/2)
  }
  +
  \nonumber\\
  &   &
  -\,
  \frac{1}{4\pi}\,
  \PV\int_{-\pi}^{\pi}d\theta_{1}\,
  \frac
  {
    \sin\!\left[\rule{0em}{2ex}(N+1/2)\Delta\theta\right]
    u(1,\theta_{1})
    -
    \cos\!\left[\rule{0em}{2ex}(N+1/2)\Delta\theta\right]
    v(1,\theta_{1})
  }
  {\sin(\Delta\theta/2)}
  +
  \nonumber\\
  &   &
  -\,
  \frac{\ii}{4\pi}\,
  \PV\int_{-\pi}^{\pi}d\theta_{1}\,
  \frac
  {
    \sin\!\left[\rule{0em}{2ex}(N+1/2)\Delta\theta\right]
    v(1,\theta_{1})
    +
    \cos\!\left[\rule{0em}{2ex}(N+1/2)\Delta\theta\right]
    u(1,\theta_{1})
  }
  {\sin(\Delta\theta/2)}.
\end{eqnarray}

\noindent
As one can see in this equation, this time all the terms involving
$\sin\!\left[\rule{0em}{2ex}(N+1/2)\Delta\theta\right]$ chancel off, and
therefore we are left with

\noindent
\begin{eqnarray}
  R_{N}(1,\theta)
  & = &
  \frac{1}{2\pi}\,
  \PV\int_{-\pi}^{\pi}d\theta_{1}\,
  \frac
  {
    \cos\!\left[\rule{0em}{2ex}(N+1/2)\Delta\theta\right]
    v(1,\theta_{1})
  }
  {
    \sin(\Delta\theta/2)
  }
  +
  \nonumber\\
  &   &
  -\,
  \frac{\ii}{2\pi}\,
  \PV\int_{-\pi}^{\pi}d\theta_{1}\,
  \frac
  {
    \cos\!\left[\rule{0em}{2ex}(N+1/2)\Delta\theta\right]
    u(1,\theta_{1})
  }
  {
    \sin(\Delta\theta/2)
  },
\end{eqnarray}

\noindent
where $\Delta\theta=\theta_{1}-\theta$ and $N\in\{1,2,3,\ldots,\infty\}$.
We have therefore our results for the real remainders, for both
$u(1,\theta)$ and $v(1,\theta)$,

\noindent
\begin{eqnarray}
  R_{N}^{F,u}(\theta)
  & = &
  \Re[R_{N}(1,\theta)]
  \nonumber\\
  & = &
  \frac{1}{2\pi}\,
  \PV\int_{-\pi}^{\pi}d\theta_{1}\,
  \frac
  {
    \cos\!\left[\rule{0em}{2ex}(N+1/2)\Delta\theta\right]
  }
  {
    \sin(\Delta\theta/2)
  }\,
  v(1,\theta_{1}),
  \nonumber\\
  R_{N}^{F,v}(\theta)
  & = &
  \Im[R_{N}(1,\theta)]
  \nonumber\\
  & = &
  -\,
  \frac{1}{2\pi}\,
  \PV\int_{-\pi}^{\pi}d\theta_{1}\,
  \frac
  {
    \cos\!\left[\rule{0em}{2ex}(N+1/2)\Delta\theta\right]
  }
  {
    \sin(\Delta\theta/2)
  }\,
  u(1,\theta_{1}),
\end{eqnarray}

\noindent
where $\Delta\theta=\theta_{1}-\theta$ and $N\in\{1,2,3,\ldots,\infty\}$.
Recalling that $f(\theta)=u(1,\theta)$ and that $g(\theta)=v(1,\theta)$
almost everywhere over the unit circle, this completes the proof of
Theorem~\ref{Theo05}.

\vspace{2.6ex}

\noindent
We believe that these are new results, written in terms of integrals which
are similar to the Dirichlet integrals, but not identical to them. Note
that the remainder of the series of $f(\theta)$ is given as an integral
involving its Fourier-conjugate function $g(\theta)$, and vice versa.
Therefore, we conclude that the convergence condition of the Fourier
series of a given real function does not depend directly on that function,
but only indirectly, through the properties of its Fourier-conjugate real
function.

Since we know that these two real functions are related by the compact
Hilbert transform, we may write these equations as

\noindent
\begin{eqnarray}
  R_{N}^{F,u}(\theta)
  & = &
  \frac{1}{2\pi}\,
  \PV\int_{-\pi}^{\pi}d\theta_{1}\,
  \frac
  {
    \cos\!\left[\rule{0em}{2ex}(N+1/2)\Delta\theta\right]
  }
  {
    \sin(\Delta\theta/2)
  }\,
  \HH[u(1,\theta_{1})],
  \nonumber\\
  R_{N}^{F,v}(\theta)
  & = &
  \frac{1}{2\pi}\,
  \PV\int_{-\pi}^{\pi}d\theta_{1}\,
  \frac
  {
    \cos\!\left[\rule{0em}{2ex}(N+1/2)\Delta\theta\right]
  }
  {
    \sin(\Delta\theta/2)
  }\,
  \HH[v(1,\theta_{1})],
\end{eqnarray}

\noindent
where $\Delta\theta=\theta_{1}-\theta$ and $N\in\{1,2,3,\ldots,\infty\}$.
Note that the two results are now identical in form. Therefore, given an
arbitrary zero-average integrable real function $f(\theta)$ on the unit
circle, we have our final result for the remainder of its Fourier series,

\begin{equation}
  R_{N}^{F}(\theta)
  =
  \frac{1}{2\pi}\,
  \PV\int_{-\pi}^{\pi}d\theta_{1}\,
  \frac
  {
    \cos\!\left[\rule{0em}{2ex}(N+1/2)\Delta\theta\right]
  }
  {
    \sin(\Delta\theta/2)
  }\,
  \HH[f(\theta_{1})],
\end{equation}

\noindent
where $\Delta\theta=\theta_{1}-\theta$ and $N\in\{1,2,3,\ldots,\infty\}$.

Once again, it is interesting to observe that this relation can be
interpreted as a linear integral operator acting on the space of
integrable zero-average real functions defined on the unit circle. The
operator $\DC[N,g(\theta)]$, acting on the compact Hilbert transform
$g(\theta)$ of such a function, results in the $N^{\rm th}$ remainder of
the Fourier series of the original function $f(\theta)$, a remainder
which, if it exists at all, is itself a zero-average integrable real
function. The integration kernel of the integral operator depends only on
$N$ and on the difference $\theta-\theta_{1}$, and is given by

\begin{equation}\label{EQIntKernDc}
  K_{\DC}(N,\theta-\theta_{1})
  =
  \frac{1}{2\pi}\,
  \frac
  {\cos\!\left[\rule{0em}{2ex}(N+1/2)(\theta_{1}-\theta)\right]}
  {\sin\!\left[\rule{0em}{2ex}(\theta_{1}-\theta)/2\right]},
\end{equation}

\noindent
where $N\in\{1,2,3,\ldots,\infty\}$, so that the action of the operator on
an arbitrarily given zero-average integrable real function $g(\theta)$ can
be written as

\begin{equation}
  \DC[N,g(\theta)]
  =
  \PV\int_{-\pi}^{\pi}d\theta_{1}\,
  K_{\DC}(N,\theta-\theta_{1})
  g(\theta_{1}).
\end{equation}

\noindent
This new operator, which we might refer to as the {\em conjugate Dirichlet
  operator}, is similar to the Dirichlet operator, and is such that the
$N^{\rm th}$ remainder of the Fourier series of the real function
$f(\theta)$ is given by the action of this operator on the
Fourier-conjugate function $g(\theta)$ of the real function $f(\theta)$,

\begin{equation}
  R_{N}^{F}(\theta)
  =
  \DC[N,g(\theta)],
\end{equation}

\noindent
where $N\in\{1,2,3,\ldots,\infty\}$. Note once more that
$\DC[N,g(\theta)]$ constitutes in fact a whole collection of linear
integral operators acting on the space of zero-average integrable real
functions. Note also that, since by hypothesis $f(\theta)$ and $g(\theta)$
are integrable on the unit circle, the Cauchy principal value refers only
to the explicit non-integrable singularity of the integration kernel at
$\theta_{1}=\theta$. In this operator notation we have therefore that the
remainder of the Fourier series of an arbitrarily given zero-average
integrable real function $f(\theta)$ is given by the composition of
$\DC[N,g(\theta)]$ with $\HH[f(\theta)]$,

\begin{equation}
  R_{N}^{F}(\theta)
  =
  \DC\!\left[\rule{0em}{2ex}N,\HH[f(\theta)]\right].
\end{equation}

\noindent
Note that, according to the inverse of the relation shown in
Equation~(\ref{EQRemMap}) we may write as well that

\begin{equation}
  R_{N}^{F}(\theta)
  =
  \HH^{-1}\!\left[\rule{0em}{2ex}\DC[N,f(\theta)]\right],
\end{equation}

\noindent
which constitutes an equivalent way to express the remainder of the
Fourier series of $f(\theta)$ in terms of the function itself.

Finally note that, given the results obtained here for the partial sums
and remainders of the Fourier series, the expressions in the infinite
collection of identities shown in Equation~(\ref{EQInfColIdenII}) have
now, in fact, the very simple interpretation that was alluded to there,
since we now see that they can in fact be written as

\noindent
\begin{eqnarray}
  f(\theta)
  & = &
  S_{N}^{F,f}(\theta)
  +
  R_{N}^{F,f}(\theta),
  \nonumber\\
  g(\theta)
  & = &
  S_{N}^{F,g}(\theta)
  +
  R_{N}^{F,g}(\theta),
\end{eqnarray}

\noindent
where $N\in\{1,2,3,\ldots,\infty\}$, a fact which greatly clarifies the
nature of that infinite collections of identities.

\section{The Convergence Condition}\label{Sec07}

Given an arbitrary zero-average integrable real function $f(\theta)$
defined on the unit circle, the necessary and sufficient condition for the
convergence of its Fourier series at the point $\theta$ is stated very
simply as the condition that

\begin{equation}
  \lim_{N\to\infty}
  R_{N}^{F}(\theta)
  =
  0,
\end{equation}

\noindent
where $R_{N}^{F}(\theta)$ is the remainder of that Fourier series, as
given in Section~\ref{Sec06}. According to what was shown in that Section,
in terms of the integral operator $\DC[N,g(\theta)]$ this translates
therefore as the condition that

\begin{equation}
  \lim_{N\to\infty}
  \DC[N,g(\theta)]
  =
  0,
\end{equation}

\noindent
where $g(\theta)$ is the Fourier-conjugate function of $f(\theta)$, which
is given by the compact Hilbert transform $\HH[f(\theta)]$, leading
therefore to the composition of the two operators,

\begin{equation}
  \lim_{N\to\infty}
  \DC\!\left[\rule{0em}{2ex}N,\HH[f(\theta)]\right]
  =
  0.
\end{equation}

\noindent
Equivalently, we may define the linear integral operator
$\DR[N,f(\theta)]$ to be this composition of $\DC[N,f(\theta)]$ with
$\HH[f(\theta)]$,

\begin{equation}
  \DR[N,f(\theta)]
  =
  \DC\!\left[\rule{0em}{2ex}N,\HH[f(\theta)]\right],
\end{equation}

\noindent
that therefore maps $f(\theta)$ directly onto the $N^{\rm th}$ remainder
of its Fourier series,

\begin{equation}
  R_{N}^{F}(\theta)
  =
  \DR[N,f(\theta)],
\end{equation}

\noindent
so that the convergence condition of the Fourier series of $f(\theta)$ can
now be written as

\begin{equation}
  \lim_{N\to\infty}
  \DR[N,f(\theta)]
  =
  0.
\end{equation}

\noindent
Combining the integration kernels of the operators $\DC[N,f(\theta)]$,
given in Equation~(\ref{EQIntKernDc}), and $\HH[f(\theta)]$, given in
Equation~(\ref{EQIntKernHc}), we may write an integration kernel for the
operator $\DR[N,f(\theta)]$, which is not given explicitly, but rather
remains expressed as an integral over the unit circle,

\begin{equation}\label{EQRemInt1}
  K_{\DR}(N,\theta,\theta_{1})
  =
  \frac{1}{4\pi^{2}}\,
  \PV\int_{-\pi}^{\pi}d\theta'\,
  \frac
  {
    \cos\!\left[\rule{0em}{2ex}(N+1/2)(\theta'-\theta)\right]
    \cos\!\left[\rule{0em}{2ex}(\theta'-\theta_{1})/2\right]
  }
  {
    \sin\!\left[\rule{0em}{2ex}(\theta'-\theta)/2\right]
    \sin\!\left[\rule{0em}{2ex}(\theta'-\theta_{1})/2\right]
  },
\end{equation}

\noindent
in terms of which the action of the operator $\DR[N,f(\theta)]$ on
$f(\theta)$ is given by

\begin{equation}
  \DR[N,f(\theta)]
  =
  \PV\int_{-\pi}^{\pi}d\theta_{1}\,
  K_{\DR}(N,\theta,\theta_{1})
  f(\theta_{1}).
\end{equation}

\noindent
Note that at this point it is not clear whether or not the integration
kernel depends only on the difference $\theta-\theta_{1}$. We will prove
that it does, and we will also write it in a somewhat more convenient
form. In this section we will prove the following theorem.

\begin{theorem}\Colon\label{Theo06}
  Given an arbitrary zero-average integrable real function $f(\theta)$
  defined on the unit circle, the necessary and sufficient condition for
  the convergence of its Fourier series at the point $\theta$ is as
  follows:
\end{theorem}

\noindent
\begin{equation}\label{EQNecSufCond}
  \lim_{N\to\infty}
  \PV\int_{-\pi}^{\pi}d\theta_{1}\,
  K_{\DR}(N,\theta-\theta_{1})
  f(\theta_{1})
  =
  0,
\end{equation}

\noindent
where the integration kernel is given by

\noindent
\begin{eqnarray}\label{EQRemInt2}
  \lefteqn
  {
    K_{\DR}(N,\theta-\theta_{1})
  }
  &   &
  \nonumber\\
  & = &
  \frac
  {\cos\!\left[\rule{0em}{2ex}(N+1)(\theta-\theta_{1})/2\right]}
  {4\pi^{2}}\,
  \PV\int_{-\pi}^{\pi}d\theta'\,
  \frac
  {
    \cos(N\theta')
  }
  {
    \cos\!\left[\rule{0em}{2ex}(\theta-\theta_{1})/2\right]
    -
    \cos(\theta')
  }
  +
  \nonumber\\
  &   &
  \hspace{1.75em}
  +
  \frac
  {\cos\!\left[\rule{0em}{2ex}N(\theta-\theta_{1})/2\right]}
  {4\pi^{2}}\,
  \PV\int_{-\pi}^{\pi}d\theta'\,
  \frac
  {
    \cos\!\left[\rule{0em}{2ex}(N+1)\theta'\right]
  }
  {
    \cos\!\left[\rule{0em}{2ex}(\theta-\theta_{1})/2\right]
    -
    \cos(\theta')
  }.
\end{eqnarray}

\vspace{2.6ex}

\noindent
Note that the two integrals in this form of the condition are almost
identical, differing only by the exchange of $N$ for $N+1$. Note also that
the integrands of these integrals are singular at the points
$\theta'=\pm(\theta-\theta_{1})/2$. Note, finally, that the condition in
Equation~(\ref{EQNecSufCond}) means that the remainder $R_{N}^{F}(\theta)$
must exist, being a finite number for each $N$, as well as that its
$N\to\infty$ limit must be zero. The existence of the remainder is, of
course, equivalent to the existence of the integrals involved.

\begin{proof}\Colon
\end{proof}

\noindent
We start by making in the integral in Equation~(\ref{EQRemInt1}) the
transformation of variables

\noindent
\begin{eqnarray}
  \theta''
  & = &
  \theta'
  -
  \frac{\theta+\theta_{1}}{2}
  \;\;\;\Rightarrow
  \nonumber\\
  \theta'
  & = &
  \theta''
  +
  \frac{\theta+\theta_{1}}{2},
\end{eqnarray}

\noindent
which implies that

\noindent
\begin{eqnarray}
  \theta'-\theta
  & = &
  \theta''
  -
  \frac{\theta-\theta_{1}}{2},
  \nonumber\\
  \theta'-\theta_{1}
  & = &
  \theta''
  +
  \frac{\theta-\theta_{1}}{2},
\end{eqnarray}

\noindent
and which also implies that $d\theta'=d\theta''$, so that we have

\noindent
\begin{eqnarray}
  \lefteqn
  {
    K_{\DR}(N,\theta,\theta_{1})
  }
  &   &
  \nonumber\\
  & = &
  \frac{1}{4\pi^{2}}\,
  \PV\int_{-\pi}^{\pi}d\theta''\,
  \frac
  {
    \cos\!\left[\rule{0em}{2ex}N_{1}\theta''-N_{1}(\theta-\theta_{1})/2\right]
    \cos\!\left[\rule{0em}{2ex}\theta''/2+(\theta-\theta_{1})/4\right]
  }
  {
    \sin\!\left[\rule{0em}{2ex}\theta''/2-(\theta-\theta_{1})/4\right]
    \sin\!\left[\rule{0em}{2ex}\theta''/2+(\theta-\theta_{1})/4\right]
  },
\end{eqnarray}

\noindent
where $N_{1}=N+1/2$ with $N\in\{1,2,3,\ldots,\infty\}$, and where we do
not have to change the integration limits since the integration runs over
a circle. Note that at this point it is already clearly apparent that
$K_{\DR}(N,\theta,\theta_{1})$ depends only on $N$ and on the difference
$\theta-\theta_{1}$, and therefore from now on we will write it as
$K_{\DR}(N,\theta-\theta_{1})$. Changing $\theta''$ back to $\theta'$ and
using the notation $\gamma=(\theta-\theta_{1})/2$ we have

\noindent
\begin{eqnarray}\label{EQKernRem}
  K_{\DR}(N,\theta-\theta_{1})
  & = &
  \frac{1}{4\pi^{2}}\,
  \PV\int_{-\pi}^{\pi}d\theta'\,
  \frac
  {
    \cos\!\left[\rule{0em}{2ex}N_{1}(\theta'-\gamma)\right]
    \cos\!\left[\rule{0em}{2ex}(\theta'+\gamma)/2\right]
  }
  {
    \sin\!\left[\rule{0em}{2ex}(\theta'-\gamma)/2\right]
    \sin\!\left[\rule{0em}{2ex}(\theta'+\gamma)/2\right]
  }
  \nonumber\\
  & = &
  \frac{1}{4\pi^{2}}\,
  \PV\int_{-\pi}^{\pi}d\theta'\,
  \frac{P(N,\theta',\gamma)}{Q(\theta',\gamma)},
\end{eqnarray}

\noindent
where $\gamma=(\theta-\theta_{1})/2$ and $N_{1}=N+1/2$ with
$N\in\{1,2,3,\ldots,\infty\}$. We will now manipulate the denominator
$Q(\theta',\gamma)$ and the numerator $P(N,\theta',\gamma)$ in this
integrand using trigonometric identities. We start with the denominator,
and using the trigonometric identities for the sum of two angles we get

\noindent
\begin{eqnarray}
  Q(\theta',\gamma)
  & = &
  \sin\!\left[\rule{0em}{2ex}(\theta'-\gamma)/2\right]
  \sin\!\left[\rule{0em}{2ex}(\theta'+\gamma)/2\right]
  \nonumber\\
  & = &
  \left[
    \rule{0em}{2ex}
    \sin(\theta'/2)
    \cos(\gamma/2)
    -
    \cos(\theta'/2)
    \sin(\gamma/2)
  \right]
  \times
  \nonumber\\
  &   &
  \times
  \left[
    \rule{0em}{2ex}
    \sin(\theta'/2)
    \cos(\gamma/2)
    +
    \cos(\theta'/2)
    \sin(\gamma/2)
  \right]
  \nonumber\\
  & = &
  \sin^{2}(\theta'/2)
  \cos^{2}(\gamma/2)
  -
  \cos^{2}(\theta'/2)
  \sin^{2}(\gamma/2).
\end{eqnarray}

\noindent
If we now write the cosines in terms of the corresponding sines we have

\noindent
\begin{eqnarray}
  Q(\theta',\gamma)
  & = &
  \sin^{2}(\theta'/2)
  -
  \sin^{2}(\theta'/2)
  \sin^{2}(\gamma/2)
  +
  \nonumber\\
  &   &
  -
  \sin^{2}(\gamma/2)
  +
  \sin^{2}(\theta'/2)
  \sin^{2}(\gamma/2)
  \nonumber\\
  & = &
  \sin^{2}(\theta'/2)
  -
  \sin^{2}(\gamma/2).
\end{eqnarray}

\noindent
Using now the half-angle trigonometric identities we finally have

\noindent
\begin{eqnarray}
  Q(\theta',\gamma)
  & = &
  \frac{1-\cos(\theta')}{2}
  -
  \frac{1-\cos(\gamma)}{2}
  \nonumber\\
  & = &
  \frac{\cos(\gamma)-\cos(\theta')}{2}.
\end{eqnarray}

\noindent
It is important to note that this is an {\em even} function of $\theta'$.
Turning now to the numerator $P(N,\theta',\gamma)$, and using the
trigonometric identities for the sum of two angles we get

\noindent
\begin{eqnarray}
  P(N,\theta',\gamma)
  & = &
  \cos\!\left[\rule{0em}{2ex}N_{1}(\theta'-\gamma)\right]
  \cos\!\left[\rule{0em}{2ex}(\theta'+\gamma)/2\right]
  \nonumber\\
  & = &
  \left[
    \rule{0em}{2ex}
    \cos(N_{1}\theta')
    \cos(N_{1}\gamma)
    +
    \sin(N_{1}\theta')
    \sin(N_{1}\gamma)
  \right]
  \times
  \nonumber\\
  &   &
  \times
  \left[
    \rule{0em}{2ex}
    \cos(\theta'/2)
    \cos(\gamma/2)
    -
    \sin(\theta'/2)
    \sin(\gamma/2)
  \right]
  \nonumber\\
  & = &
  \cos(N_{1}\theta')
  \cos(N_{1}\gamma)
  \cos(\theta'/2)
  \cos(\gamma/2)
  +
  \nonumber\\
  &   &
  -
  \cos(N_{1}\theta')
  \cos(N_{1}\gamma)
  \sin(\theta'/2)
  \sin(\gamma/2)
  +
  \nonumber\\
  &   &
  +
  \sin(N_{1}\theta')
  \sin(N_{1}\gamma)
  \cos(\theta'/2)
  \cos(\gamma/2)
  +
  \nonumber\\
  &   &
  -
  \sin(N_{1}\theta')
  \sin(N_{1}\gamma)
  \sin(\theta'/2)
  \sin(\gamma/2).
\end{eqnarray}

\noindent
Of these four terms, the first and last ones are even on $\theta'$, and
the two middle ones are odd. Since the denominator is even and the
integral on $\theta'$ shown in Equation~(\ref{EQKernRem}) is over a
symmetric interval, the integrals of the two middle terms will be zero,
and therefore we can ignore these two terms of the numerator. We thus
obtain for our kernel

\begin{equation}
  K_{\DR}(N,\theta-\theta_{1})
  =
  \frac{1}{4\pi^{2}}\,
  \PV\int_{-\pi}^{\pi}d\theta'\,
  \frac{T(N,\theta',\gamma)}{Q(\theta',\gamma)},
\end{equation}

\noindent
where the new numerator is given by

\noindent
\begin{eqnarray}
  T(N,\theta',\gamma)
  & = &
  \left\{
    \rule{0em}{2.5ex}
    \cos\!\left[\rule{0em}{2ex}(N+1/2)\theta'\right]
    \cos(\theta'/2)
  \right\}
  \cos(N_{1}\gamma)
  \cos(\gamma/2)
  +
  \nonumber\\
  &   &
  -
  \left\{
    \rule{0em}{2.5ex}
    \sin\!\left[\rule{0em}{2ex}(N+1/2)\theta'\right]
    \sin(\theta'/2)
  \right\}
  \sin(N_{1}\gamma)
  \sin(\gamma/2),
\end{eqnarray}

\noindent
where $\gamma=(\theta-\theta_{1})/2$ and $N_{1}=N+1/2$ with
$N\in\{1,2,3,\ldots,\infty\}$. Using once again the trigonometric
identities for the sum of two angles in order to work on the two
expressions within pairs of curly brackets, we get

\noindent
\begin{eqnarray}
  B_{\rm c}
  & = &
  \left\{
    \rule{0em}{2.5ex}
    \cos\!\left[\rule{0em}{2ex}(N+1/2)\theta'\right]
    \cos(\theta'/2)
  \right\}
  \nonumber\\
  & = &
  \left[
    \rule{0em}{2ex}
    \cos(N\theta')
    \cos(\theta'/2)
    -
    \sin(N\theta')
    \sin(\theta'/2)
  \right]
  \cos(\theta'/2)
  \nonumber\\
  & = &
  \rule{0em}{2ex}
  \cos(N\theta')
  \cos^{2}(\theta'/2)
  -
  \sin(N\theta')
  \sin(\theta'/2)
  \cos(\theta'/2),
  \nonumber\\
  B_{\rm s}
  & = &
  \left\{
    \rule{0em}{2.5ex}
    \sin\!\left[\rule{0em}{2ex}(N+1/2)\theta'\right]
    \sin(\theta'/2)
  \right\}
  \nonumber\\
  & = &
  \left[
    \rule{0em}{2ex}
    \sin(N\theta')
    \cos(\theta'/2)
    +
    \cos(N\theta')
    \sin(\theta'/2)
  \right]
  \sin(\theta'/2)
  \nonumber\\
  & = &
  \rule{0em}{2ex}
  \sin(N\theta')
  \cos(\theta'/2)
  \sin(\theta'/2)
  +
  \cos(N\theta')
  \sin^{2}(\theta'/2).
\end{eqnarray}

\noindent
Using now the half-angle trigonometric identities in order to eliminate
the functions $\cos(\theta'/2)$ and $\sin(\theta'/2)$ in favor of
$\cos(\theta')$ and $\sin(\theta')$, we get

\noindent
\begin{eqnarray}
  B_{\rm c}
  & = &
  \rule{0em}{2ex}
  \cos(N\theta')\,
  \frac{1+\cos(\theta')}{2}
  -
  \sin(N\theta')\,
  \frac{\sin(\theta')}{2}
  \nonumber\\
  & = &
  \frac{\cos(N\theta')}{2}
  +
  \frac
  {
    \cos(N\theta')\cos(\theta')
    -
    \sin(N\theta')\sin(\theta')
  }
  {2},
  \nonumber\\
  B_{\rm s}
  & = &
  \rule{0em}{2ex}
  \sin(N\theta')\,
  \frac{\sin(\theta')}{2}
  +
  \cos(N\theta')\,
  \frac{1-\cos(\theta')}{2}
  \nonumber\\
  & = &
  \frac{\cos(N\theta')}{2}
  -
  \frac
  {
    \cos(N\theta')\cos(\theta')
    -
    \sin(N\theta')\sin(\theta')
  }
  {2}.
\end{eqnarray}

\noindent
Using once more the trigonometric identities for the sum of two angles we
get

\noindent
\begin{eqnarray}
  B_{\rm c}
  & = &
  \frac{\cos(N\theta')}{2}
  +
  \frac{\cos\!\left[\rule{0em}{2ex}(N+1)\theta'\right]}{2},
  \nonumber\\
  B_{\rm s}
  & = &
  \frac{\cos(N\theta')}{2}
  -
  \frac{\cos\!\left[\rule{0em}{2ex}(N+1)\theta'\right]}{2}.
\end{eqnarray}

\noindent
We get therefore for our new numerator, using yet gain the trigonometric
identities for the sum of two angles,

\noindent
\begin{eqnarray}
  T(N,\theta',\gamma)
  & = &
  \frac{\cos(N\theta')}{2}
  \left[
    \rule{0em}{2ex}
    \cos(N_{1}\gamma)
    \cos(\gamma/2)
    -
    \sin(N_{1}\gamma)
    \sin(\gamma/2)
  \right]
  +
  \nonumber\\
  &   &
  +
  \frac{\cos\!\left[\rule{0em}{2ex}(N+1)\theta'\right]}{2}
  \left[
    \rule{0em}{2ex}
    \cos(N_{1}\gamma)
    \cos(\gamma/2)
    +
    \sin(N_{1}\gamma)
    \sin(\gamma/2)
  \right]
  \nonumber\\
  & = &
  \frac{\cos(N\theta')\cos\!\left[\rule{0em}{2ex}(N+1)\gamma\right]}{2}
  +
  \frac{\cos\!\left[\rule{0em}{2ex}(N+1)\theta'\right]\cos(N\gamma)}{2}.
\end{eqnarray}

\noindent
We therefore have for the kernel

\noindent
\begin{eqnarray}
  \lefteqn
  {
    K_{\DR}(N,\theta-\theta_{1})
  }
  &   &
  \nonumber\\
  & = &
  \frac{1}{4\pi^{2}}\,
  \PV\int_{-\pi}^{\pi}d\theta'\,
  \frac
  {
    \cos(N\theta')\cos\!\left[\rule{0em}{2ex}(N+1)\gamma\right]
    +
    \cos\!\left[\rule{0em}{2ex}(N+1)\theta'\right]\cos(N\gamma)
  }
  {
    \cos(\gamma)-\cos(\theta')
  },
\end{eqnarray}

\noindent
where $\gamma=(\theta-\theta_{1})/2$ and $N\in\{1,2,3,\ldots,\infty\}$.
This is exactly the form of the integration kernel of the operator
$\DR[N,f(\theta)]$ given in Equation~(\ref{EQRemInt2}), so that this
completes the proof of Theorem~\ref{Theo06}.

\vspace{2.6ex}

\noindent
Since the condition stated in Theorem~\ref{Theo06} is a necessary and
sufficient condition on the real function $f(\theta)$ for the convergence
of its Fourier series, any other such condition must be equivalent to it.
Note that this type of condition is not really what is usually meant by a
Fourier theorem. Those are just sufficient conditions for the convergence,
usually related to some fairly simple and easily identifiable
characteristics of the real functions, such as continuity,
differentiability, existence of lateral limits, or limited variation.
However, all such Fourier theorems must be related to this condition in
the sense that they must imply its validity.

\section{Conclusions and Outlook}\label{Sec08}

We have shown that the complex-analytic structure that we introduced
in~\cite{CAoRFI} can be used to discuss the issue of the convergence of
Fourier series. Using that structure we derived the known formulas for the
partial sums of a Fourier series, in terms of Dirichlet integrals. From
that same structure, and in fact as part of the same argument, we also
obtained a new result, namely formulas giving the remainders of a Fourier
series in terms of a similar but considerably more complex type of
integral, in fact a double integral.

The introduction of a modified version of the Hilbert transform, which we
named the compact Hilbert transform, had a central role to play in this
development. The main result that follows from it is the necessary and
sufficient condition for the convergence of a Fourier series, which is
expressed in Equations~(\ref{EQNecSufCond}) and~(\ref{EQRemInt2}). One
might consider whether or not the integration kernel involved in this new
type of Dirichlet integral, versions of which are shown in
Equations~(\ref{EQRemInt1}) and~(\ref{EQRemInt2}), can be cast in some
more convenient form, and possibly even calculated in close form in some
useful way. So far no meaningful results of this type have been found.

A direct calculation of this integration kernel in closed form seems to be
difficult, and possibly not overly useful, since it seems that such
calculations tend to just take us back to the rather trivial identity

\begin{equation}
  \DR[N,f(\theta)]
  =
  {\cal I}[f(\theta)]
  -
  \DS[N,f(\theta)],
\end{equation}

\noindent
where ${\cal I}[f(\theta)]$ is the identity operator, whose integration
kernel is a Dirac delta ``function'', an identity which is simply
equivalent to

\begin{equation}
  R_{N}^{F}(\theta)
  =
  f(\theta)
  -
  S_{N}^{F}(\theta).
\end{equation}

\noindent
In so far as can be currently ascertained, this identity does not provide
any constructively useful information about the convergence problem.

We already knew that he convergence problem of Fourier series relates to
the existence and nature of singularities of the corresponding inner
analytic functions on the unit circle. The convergence condition expressed
in Equations~(\ref{EQNecSufCond}) and~(\ref{EQRemInt2}) reflect this
relation, since the existence of non-integrable singularities at the unit
circle may very well disturb the validity of the condition by making the
integrals diverge or at least not go to zero in the $N\to\infty$
limit. Note that this relation is non-local because, due to the fact that
the integrals are over the whole unit circle, the existence of a
non-integrable singularity at a single point may prevent the convergence
of the series at almost all points.

It seems to us, at this time, that the most promising possible development
of the convergence analysis presented here is probably one targeted at
detailing the relation between the convergence of the series and the
specific classification of the singularities on the unit circle, involving
the concepts of hard and soft singularities, as well as the corresponding
degrees of hardness or softness that can be attributed to them.

\section*{Acknowledgments}

The author would like to thank his friend and colleague Prof. Carlos
Eugênio Imbassay Carneiro, to whom he is deeply indebted for all his
interest and help, as well as his careful reading of the manuscript and
helpful criticism regarding this work.

The author would also like to thank the anonymous referee of the first
paper in this series~\cite{CAoRFI}, whose objection regarding the relation
between the complex-analytic structure and the known Fourier convergence
theorems put him on a thought path that ultimately lead him to the main
results presented in this paper.

\bibliography{allrefs_en}\bibliographystyle{ieeetr}

\end{document}